\pdfoutput=1
\documentclass[11pt]{amsart}

\usepackage{amsmath,amssymb,amsthm,mathrsfs}
\usepackage{geometry}
\usepackage{esint}
\usepackage{hyperref}
\usepackage{enumitem}

\def\XXint#1#2#3{{\setbox0=\hbox{$#1{#2#3}{\int}$}
     \vcenter{\hbox{$#2#3$}}\kern-.5\wd0}}

\newtheorem{theorem}{Theorem}[section]
\newtheorem{lemma}[theorem]{Lemma}

\newtheorem{definition}[theorem]{Definition}
\newtheorem{corollary}[theorem]{Corollary}
\newtheorem{proposition}[theorem]{Proposition}
\theoremstyle{definition}

\theoremstyle{remark}
\newtheorem*{remark}{Remark}

\title{TRANSLATORS ASYMPTOTIC TO VERTICAL HYPERPLANES}

\author{Yunfan Zhao}
\email{yfzucla21@ucla.edu}

\date{} 

\begin{document}

\begin{abstract}
We prove that a vertical hyperplane is the only complete translator, properly immersed in $\mathbb{R}^{n+1}$ ($n \ge 3$) and having finite topology, whose ends are $C^1$-asymptotic to vertical hyperplanes.
\end{abstract}

\maketitle

\section{Introduction}

\subsection{Background and Literature Review}

The systematic investigation of the mean curvature flow necessitates a rigorous analysis of its singularity models, prominently including translating solitons. A smooth hypersurface $\Sigma^{n} \hookrightarrow \mathbb{R}^{n+1}$ is defined to be a translator if its mean curvature vector $H$ satisfies the structural equation
\[
H = e_{n+1}^{\perp}.
\]
These self-similar geometric entities are of paramount importance to the broader theory, as they model Type-II singularities of the mean curvature flow. The foundational analytical apparatus governing these flows, including the formulation of weak solutions via elliptic regularization, was systematically developed by Ilmanen~\cite{Ilmanen1994}. The intrinsic geometric constraints of the flow were previously elucidated through the Harnack estimates established by Hamilton~\cite{Hamilton1995} and the convexity estimates formulated by Huisken and Sinestrari~\cite{HuiskenSinestrari1999}.

Within the taxonomy of these structures, the complete classification of convex translating solutions was achieved by Wang~\cite{Wang2011}. The dynamical stability of such translating solutions under the mean curvature flow was subsequently proven by Clutterbuck, Schn\"urer, and Schulze~\cite{ClutterbuckSchnurerSchulze2007}. Further structural constraints on translating solitons have been progressively formalized: Kunikawa~\cite{Kunikawa2015} established Bernstein-type rigidity theorems for translators, while Xin~\cite{Xin2015} investigated their intrinsic geometry. The topological restrictions imposed upon translating solitons were analyzed by Mart\'in, Savas-Halilaj, and Smoczyk~\cite{MartinSavasHalilajSmoczyk2015}.

In more recent developments, analytical attention has converged on the asymptotic geometry of translator ends. The behavior of translating solitons asymptotic to hyperplanes in $\mathbb{R}^{n+1}$ was studied by Gama and Mart\'in~\cite{GamaMartin2020}. In the case of complete embedded translators, Khan~\cite{Khan2023} showed that finite total curvature forces the ends of the surface to be asymptotic to vertical planes.

Motivated by these structural results, Lynch and Tinaglia~\cite{LynchTinaglia2026} established a Liouville-type rigidity theorem: any complete translator properly immersed in $\mathbb{R}^{3}$ with finite topology, whose ends are asymptotic to vertical planes, must itself be a vertical plane. Their argument relies on a flux identity expressing the total mean curvature in terms of the boundary flux of the translation vector field.

In the present work, we extend this classification paradigm to arbitrary dimensions $n \geq 3$. The global rigidity is established through a precise analysis of the asymptotic decay of solutions, obtained via an algebraic decomposition of the quasilinear translator equation and a functional analysis of the associated geometric flux.

\subsection{Main Contributions}

Our main theorem classifies complete, properly immersed, connected translators with finite topology whose ends are $C^1$-asymptotic to vertical hyperplanes in ambient dimensions $n \ge 3$.

\medskip

\noindent \textbf{Theorem~\ref{thm:Rigidity of Asymptotic Translators}}
Let $\Sigma^n \hookrightarrow \mathbb{R}^{n+1}$ ($n \ge 3$) be a complete, properly immersed, connected translator with finite topology whose ends are $C^1$-asymptotic to vertical hyperplanes. Then $\Sigma$ is globally a vertical hyperplane.

\subsection{Method of Approach}

Our proof of Theorem~\ref{thm:Rigidity of Asymptotic Translators} relies principally on an identity expressing the total mean curvature of any admissible translator domain as the flux of the translation vector $e_{n+1}^\top$ across its boundary. To demonstrate that this boundary flux vanishes in the limit along a cylindrical exhaustion of the ends, we introduce a strict measurable partition of the asymptotic exterior domain into a canonical wake set and an exponential region. This partition allows us to establish bifurcated, pointwise decay estimates for the gradient using the drift Laplacian $L_0 = \Delta + \partial_{y_n}$. The vanishing boundary flux implies that the mean curvature identically vanishes, reducing $\Sigma$ to an asymptotic minimal graph. We complete the proof by proving a Liouville rigidity theorem for such minimal graphs, which strictly forces the solution to be flat.

\section{Preliminaries}

The systematic study of translator hypersurfaces necessitates, as a prerequisite, the establishment of a rigorous formal language and the precise delineation of the ambient geometric spaces. This chapter furnishes the requisite analytic and geometric apparatus. We construct the foundational definitions governing structural translator equations and asymptotic behaviors. Furthermore, we introduce the precise spatial stratifications and specialized functional spaces that will subsequently permit the systematic domination and iterative absorption of nonlinear error terms in our asymptotic analysis.

\subsection{Basic Definitions}

Before undertaking any global analysis, one must strictly define the primary objects of study. Here, we formalize the notion of a translator as a smooth hypersurface $\Sigma^n \hookrightarrow \mathbb{R}^{n+1}$ satisfying the fundamental structural equation $H = e_{n+1}^{\perp}$. Additionally, we rigorously define the asymptotic limit of its ends, specifying the exact $C^1$ decay conditions—namely, the vanishing of the gradient norm $\|Du(y)\|$—required for an end to be considered asymptotic to a canonically normalized vertical hyperplane. This establishes the initial topological and metric boundary conditions upon which our eventual classification depends.

\begin{definition}
A smooth hypersurface $\Sigma^n \hookrightarrow \mathbb{R}^{n+1}$ is called a \textit{translator} if its mean curvature vector $\mathbf{H}$ satisfies the structural equation:
\begin{equation}
\mathbf{H} = e_{n+1}^\perp
\end{equation}
where $(\cdot)^\perp$ denotes the projection onto the normal bundle $N\Sigma$.
\end{definition}

\begin{definition}[Vertical Hyperplane]
A hyperplane $\Pi \subset \mathbb{R}^{n+1}$ is termed \textit{vertical} if it contains the translation vector, i.e., $e_{n+1} \in \Pi$. Up to rigid motion, we normalize $\Pi$ to be the coordinate plane $\{x_1 = 0\}$. We adopt coordinates $y = (y_1, \dots, y_n)$ on $\Pi$, identifying the $y_n$-axis with the $x_{n+1}$-direction.
\end{definition}

\begin{definition}[Asymptotic End]\label{def:Asymptotic End}
An end $E \subset \Sigma$ is said to be \textit{$C^1$-asymptotic to a vertical hyperplane $\Pi$} if there exists a compact set $K \subset \Pi$ and a function $u \in C^\infty(\Pi \setminus K)$ such that $E = \text{graph}(u)$ and the gradient satisfies the decay condition:
\begin{equation}
\lim_{|y| \to \infty} \|Du(y)\| = 0
\end{equation}
\end{definition}

\subsection{Spatial Partitions}

To govern the asymptotic behavior of solutions, it is indispensable to dissect the underlying Euclidean exterior domain into distinct regimes dictated by the nature of the ambient drift. We define a strict Lebesgue-measurable partition of the spatial hyperplane, bifurcating the space into a Canonical Wake Set $\mathcal{W}_C$ and a Canonical Exponential Set $\mathcal{E}_C$. This algebraic and measurable stratification is of paramount significance, as it strictly dictates the bifurcated, pointwise decay estimates for the gradient during subsequent functional analysis.

\begin{definition}(Strict Asymptotic Boundedness)
Let $E$ be a topological vector space and $\Omega \subset \mathbb{R}^n$ be an exterior domain. A function $f: \Omega \rightarrow E$ is said to be strictly asymptotically bounded by a non-negative function $g: \Omega \rightarrow \mathbb{R}_{\geq 0}$, denoted $f \in \mathcal{B}_{\infty}(g)$, if there exists a strictly positive real constant $C>0$ and a compact set $K \subset \mathbb{R}^n$ such that for all $y \in \Omega \backslash K$, the norm satisfies $\|f(y)\|_E \leq C g(y)$.
\end{definition}

\begin{definition}[Spatial Partition of the Hyperplane]
Let $\Pi \cong \mathbb{R}^n$ be a Euclidean hyperplane and $e_n$ be the translation unit vector. For any strictly positive real constant $C>0$, we define the Canonical Wake Set $\mathcal{W}_C \subset \mathbb{R}^n$ and the Canonical Exponential Set $\mathcal{E}_C \subset \mathbb{R}^n$ by the strict algebraic relations:

\begin{equation} 
\begin{aligned}
\mathcal{W}_C & :=\left\{y \in \mathbb{R}^n:|y|+y_n \leq C\right\} \\
\mathcal{E}_C & :=\left\{y \in \mathbb{R}^n:|y|+y_n>C\right\}
\end{aligned}
\end{equation}
\end{definition}

\begin{lemma}[Spherical Measurable Partition]\label{lem:Spherical Measurable Partition}
For any strictly positive real constant $C > 0$ and any radius $R > 0$, the spatial intersections $\mathcal{W}_{C, R} := \mathcal{W}_{C} \cap \partial B_{R}(0)$ and $\mathcal{E}_{C, R} := \mathcal{E}_{C} \cap \partial B_{R}(0)$ constitute a strict Lebesgue-measurable partition of the Euclidean sphere $\mathcal{S}_{R} = \partial B_{R}(0)$ with respect to the standard spherical surface measure $d\sigma$.
\end{lemma}

\begin{proof}
We define the phase mapping $g: \mathbb{R}^{n} \to \mathbb{R}$ by $g(y) = |y| + y_{n}$. By the definitions of the Canonical sets, $\mathcal{W}_{C} = g^{-1}((-\infty, C])$ and $\mathcal{E}_{C} = g^{-1}((C, \infty))$.

The intervals $I_{1} = (-\infty, C]$ and $I_{2} = (C, \infty)$ form a strict algebraic partition of the codomain $\mathbb{R}$ (i.e., $I_{1} \cup I_{2} = \mathbb{R}$ and $I_{1} \cap I_{2} = \emptyset$). Because the preimage operation $g^{-1}$ preserves Boolean set operations, we  deduce:
\begin{equation}
\mathcal{W}_{C} \cup \mathcal{E}_{C} = g^{-1}(I_{1} \cup I_{2}) = g^{-1}(\mathbb{R}) = \mathbb{R}^{n}
\end{equation}
\begin{equation}
\mathcal{W}_{C} \cap \mathcal{E}_{C} = g^{-1}(I_{1} \cap I_{2}) = g^{-1}(\emptyset) = \emptyset
\end{equation}
Thus, $\mathcal{W}_{C}$ and $\mathcal{E}_{C}$ form a strict algebraic partition of $\mathbb{R}^{n}$.

$g$ is a continuous mapping. In the standard topology of $\mathbb{R}$, the ray $(-\infty, C]$ is closed and the ray $(C, \infty)$ is open. Because the continuous preimage of a closed set is closed, $\mathcal{W}_{C}$ is topologically closed in $\mathbb{R}^{n}$. Because the continuous preimage of an open set is open, $\mathcal{E}_{C}$ is topologically open in $\mathbb{R}^{n}$.

The Borel $\sigma$-algebra $\mathcal{B}(\mathbb{R}^{n})$ is generated by the open sets, and is closed under complementation (thus containing all closed sets). Therefore, we strictly conclude $\mathcal{W}_{C} \in \mathcal{B}(\mathbb{R}^{n})$ and $\mathcal{E}_{C} \in \mathcal{B}(\mathbb{R}^{n})$.

We restrict our domain to the subspace $\mathcal{S}_{R} = \partial B_{R}(0) \subset \mathbb{R}^{n}$. By Boolean intersection distribution on the global partition established:
\begin{equation}
(\mathcal{W}_{C} \cap \mathcal{S}_{R}) \cup (\mathcal{E}_{C} \cap \mathcal{S}_{R}) = (\mathcal{W}_{C} \cup \mathcal{E}_{C}) \cap \mathcal{S}_{R} = \mathbb{R}^{n} \cap \mathcal{S}_{R} = \mathcal{S}_{R}
\end{equation}
\begin{equation}
(\mathcal{W}_{C} \cap \mathcal{S}_{R}) \cap (\mathcal{E}_{C} \cap \mathcal{S}_{R}) = (\mathcal{W}_{C} \cap \mathcal{E}_{C}) \cap \mathcal{S}_{R} = \emptyset \cap \mathcal{S}_{R} = \emptyset
\end{equation}
Thus, $\mathcal{W}_{C, R}$ and $\mathcal{E}_{C, R}$ geometrically partition $\mathcal{S}_{R}$.

Finally, the standard spherical surface measure $d\sigma$ is a Borel measure strictly defined on the measurable space $(\mathcal{S}_{R}, \mathcal{B}(\mathcal{S}_{R}))$. Because $\mathcal{W}_{C} \in \mathcal{B}(\mathbb{R}^{n})$ and $\mathcal{E}_{C} \in \mathcal{B}(\mathbb{R}^{n})$, their intersections with $\mathcal{S}_{R}$ satisfy:
\begin{equation}
\mathcal{W}_{C, R} \in \mathcal{B}(\mathcal{S}_{R}) \quad \text{and} \quad \mathcal{E}_{C, R} \in \mathcal{B}(\mathcal{S}_{R})
\end{equation}
Because both sets belong to the trace $\sigma$-algebra upon which $d\sigma$ is defined, they are strictly $d\sigma$-measurable. Therefore, they constitute a measurable partition of the sphere.

\end{proof}

\subsection{Special Functions and Asymptotic Bounds}

The investigation of the linearized drift operator intrinsically relies upon the properties of specific transcendental kernels. In this section, we introduce the formal integral construction of the modified Bessel function of the second kind, denoted $K_{\nu}(z)$. We establish strict, a priori asymptotic bounds for this function. The analytical significance of these rigorous estimates lies in their direct application to bounding the fundamental solution of the drift equation, which is structurally indispensable for constructing representation formulas later.

\begin{definition}(Integral Construction of the Bessel Function)
Let $\nu \in \mathbb{R}$. For any strictly positive real number $z \in \mathbb{R}_{>0}$, the modified Bessel function of the second kind, denoted $K_\nu(z)$, is defined by the absolutely convergent integral:

\begin{equation} 
K_\nu(z)=\int_0^{\infty} e^{-z \cosh t} \cosh (\nu t) d t
\end{equation}
\end{definition}

\begin{proposition}[Strict Asymptotic Bound of the Bessel Function]
Let $\nu \in \mathbb{R}$. For $z>0$ define the modified Bessel function of the second kind by the absolutely convergent integral
\begin{equation}
K_\nu(z)=\int_0^\infty e^{-z\cosh t}\cosh(\nu t)\,dt .
\end{equation}
Then there exist constants $C_0>0$ and $z_0>0$ such that for all $z>z_0$,
\begin{equation}
\left|K_\nu(z)-\sqrt{\frac{\pi}{2z}}e^{-z}\right|\le C_0 z^{-3/2}e^{-z}.
\end{equation}
In particular there exists $C_1>0$ such that $K_\nu(z)\le C_1 z^{-1/2}e^{-z}$ for all sufficiently large $z$.
\end{proposition}

\begin{proof}
Let $s=\cosh t-1$. Then $t\in(0,\infty)$ maps bijectively to $s\in(0,\infty)$. The differential satisfies $ds=\sinh t\,dt$, and since $\sinh t=\sqrt{\cosh^2 t-1}=\sqrt{(s+1)^2-1}=\sqrt{s^2+2s}$, we have $dt=ds/\sqrt{2s(1+s/2)}$.

Substituting into the definition of $K_\nu$ yields
\begin{equation}
K_\nu(z)=\frac{e^{-z}}{\sqrt{2}}\int_0^\infty e^{-zs}s^{-1/2}\frac{\cosh(\nu\,\mathrm{arccosh}(1+s))}{\sqrt{1+s/2}}\,ds .
\end{equation}

Define the auxiliary function $H(s)=\frac{\cosh(\nu\,\mathrm{arccosh}(1+s))}{\sqrt{1+s/2}}$. Evaluating at $s=0$ gives $H(0)=1$. By Taylor's theorem with remainder we may write $H(s)=1+sR(s)$, where $R$ is continuous. Because $\cosh$ and $\mathrm{arccosh}$ are smooth and grow at most exponentially, there exist constants $M>0$ and $\alpha>0$ such that $|R(s)|\le M e^{\alpha s}$ for all $s>0$.

Substituting $H(s)=1+sR(s)$ into the integral gives
\begin{equation}
K_\nu(z)=\frac{e^{-z}}{\sqrt{2}}\int_0^\infty e^{-zs}s^{-1/2}\,ds
+\frac{e^{-z}}{\sqrt{2}}\int_0^\infty e^{-zs}s^{1/2}R(s)\,ds .
\end{equation}

For the first integral, make the substitution $w=zs$. Then $ds=dw/z$ and
\begin{equation}
\int_0^\infty e^{-zs}s^{-1/2}\,ds
=z^{-1/2}\int_0^\infty e^{-w}w^{-1/2}\,dw
=z^{-1/2}\Gamma(1/2)=\sqrt{\frac{\pi}{z}} .
\end{equation}
Hence the principal contribution equals $\sqrt{\frac{\pi}{2z}}e^{-z}$.

Let $E(z)$ denote the second integral. Using the bound $|R(s)|\le M e^{\alpha s}$ and assuming $z>\alpha$, we obtain
\begin{equation}
|E(z)|\le \frac{e^{-z}}{\sqrt{2}}M\int_0^\infty e^{-(z-\alpha)s}s^{1/2}\,ds .
\end{equation}

Set $u=(z-\alpha)s$. Then $ds=du/(z-\alpha)$ and
\begin{equation}
\int_0^\infty e^{-(z-\alpha)s}s^{1/2}\,ds
=(z-\alpha)^{-3/2}\int_0^\infty e^{-u}u^{1/2}\,du
=(z-\alpha)^{-3/2}\Gamma(3/2).
\end{equation}

Therefore
\begin{equation}
|E(z)|\le \frac{M\Gamma(3/2)}{\sqrt{2}}e^{-z}(z-\alpha)^{-3/2}.
\end{equation}
For $z>2\alpha$ one has $(z-\alpha)^{-3/2}\le (z/2)^{-3/2}=2^{3/2}z^{-3/2}$. Hence there exists a constant $C_0>0$ such that
\begin{equation}
|E(z)|\le C_0 z^{-3/2}e^{-z}.
\end{equation}

Combining the principal term and the error term yields
\begin{equation}
\left|K_\nu(z)-\sqrt{\frac{\pi}{2z}}e^{-z}\right|\le C_0 z^{-3/2}e^{-z}.
\end{equation}

Finally, the triangle inequality implies that for sufficiently large $z$ there exists $C_1>0$ such that
\begin{equation}
K_\nu(z)\le C_1 z^{-1/2}e^{-z}.
\end{equation}
This completes the proof.
\end{proof}

\begin{remark}
It is proper to observe that the principal asymptotic behavior of the integral, which is characterized by the exponential decay $e^{-z}$ and the inverse square-root polynomial weight $z^{-1/2}$, manifests strictly independently of the parameter $\nu$. The difference between the exact function and its principal asymptotic term is rigorously bounded by an error of order $\mathcal{O}(z^{-3/2}e^{-z})$ for sufficiently large $z$.
\end{remark}

\subsection{Weighted Norms and Iterative Analysis}

The systematic suppression of nonlinear remainders requires an appropriate topological framework for our function spaces. We define strictly positive continuous weight functions to construct both standard and exponentially weighted supremum norms over exterior domains. By formulating an iterative absorption lemma, we provide a rigorous functional-analytic mechanism to bound non-increasing functions. This abstract machinery is critical for absorbing the nonlinear functional components of the quasilinear equation into the linear principal part, thereby enabling global supremum estimates.

\begin{definition}(Weighted Supremum Norm)
Let $\Omega \subset \mathbb{R}^n$ be an exterior domain and let $w: \Omega \rightarrow \mathbb{R}_{>0}$ be a strictly positive continuous weight function. For any function $v \in L_{l o c}^{\infty}(\Omega)$, we define the weighted supremum norm over a subset $U \subset \Omega$ as:

\begin{equation} 
\|v\|_{w, U}:=\sup _{y \in U}(w(y)|v(y)|)
\end{equation}
\end{definition}

\begin{lemma}[Iterative Absorption]\label{lem:Iterative Absorption}
Let $R_{0} > 0$. Let $f:[R_{0},\infty) \rightarrow \mathbb{R}_{\ge 0}$ be a non-increasing function. Suppose there exist constants $\theta \in (0,1)$ and $C > 0$, and a strictly positive function $E:[R_{0},\infty) \rightarrow \mathbb{R}_{>0}$ satisfying the following properties.

\begin{enumerate}
    \item There exists a constant $M \ge 1$ such that for all $R \ge 2R_{0}$ we have $E(R/2) \le M E(R)$.
    \item The constants satisfy $\theta M < 1$.
    \item The function $E$ is locally bounded below on $[R_{0},2R_{0}]$ in the sense that $\inf_{s \in [R_{0},2R_{0}]} E(s) = E_{\min} > 0$.
\end{enumerate}

Assume that for all $R \ge 2R_{0}$ the inequality
\begin{equation} 
f(R) \le \theta f(R/2) + C E(R)
\end{equation}
holds. Then there exists a constant $C' > 0$, depending only on $\theta$, $C$, $M$, $E_{\min}$, and $f(R_{0})$, such that for all $R \ge R_{0}$ we have
\begin{equation} 
f(R) \le C' E(R).
\end{equation}
\end{lemma}

\begin{proof}
Let $R \ge R_{0}$. If $R_{0} \le R < 2R_{0}$, the monotonicity of $f$ gives $f(R) \le f(R_{0})$. Since $E(R) \ge E_{\min} > 0$ on this interval, we obtain $f(R) \le \frac{f(R_{0})}{E_{\min}} E(R)$.

Assume now that $R \ge 2R_{0}$. By the Archimedean property of the real numbers there exists an integer $k \ge 1$ such that $2^{k}R_{0} \le R < 2^{k+1}R_{0}$, which is equivalent to $R_{0} \le R/2^{k} < 2R_{0}$.

Iterating the inequality $f(R) \le \theta f(R/2) + C E(R)$ exactly $k$ times yields
\begin{equation} 
f(R) \le \theta^{k} f(R/2^{k}) + C \sum_{j=0}^{k-1} \theta^{j} E(R/2^{j}).
\end{equation}

We first estimate the summation term. From the scale bound $E(R/2) \le M E(R)$ it follows inductively that $E(R/2^{j}) \le M^{j} E(R)$ for every integer $j \ge 0$. Therefore
\begin{equation} 
\sum_{j=0}^{k-1} \theta^{j} E(R/2^{j}) \le E(R) \sum_{j=0}^{k-1} (\theta M)^{j}.
\end{equation}
Since $\theta M < 1$, the geometric series is bounded by $\sum_{j=0}^{k-1} (\theta M)^{j} < \sum_{j=0}^{\infty} (\theta M)^{j} = \frac{1}{1-\theta M}$. Hence the entire summation term is bounded by $\frac{C}{1-\theta M} E(R)$.

We next estimate the remaining term $\theta^{k} f(R/2^{k})$. Because $f$ is non-increasing and $R/2^{k} \ge R_{0}$, we have $f(R/2^{k}) \le f(R_{0})$. The scale bound implies $E(R/2^{k}) \le M^{k} E(R)$, which gives $(1/M)^{k} \le \frac{E(R)}{E(R/2^{k})}$. Since $R/2^{k} \in [R_{0},2R_{0})$, we also have $E(R/2^{k}) \ge E_{\min}$, and therefore $(1/M)^{k} \le \frac{E(R)}{E_{\min}}$. Because $\theta M < 1$, we have $\theta < 1/M$ and thus $\theta^{k} < (1/M)^{k}$. Combining these inequalities yields $\theta^{k} < \frac{E(R)}{E_{\min}}$, and consequently $\theta^{k} f(R/2^{k}) \le \frac{f(R_{0})}{E_{\min}} E(R)$.

Substituting the bounds for both terms into the iterated inequality gives
\begin{equation} 
f(R) \le \frac{f(R_{0})}{E_{\min}} E(R) + \frac{C}{1-\theta M} E(R).
\end{equation}
Factoring out $E(R)$ and defining $C' = \frac{f(R_{0})}{E_{\min}} + \frac{C}{1-\theta M}$ completes the proof.
\end{proof}

\begin{definition}(Exponentially Weighted Supremum Norm)
Let $\mu>0$. For a function $f$, define $\|f\|_{\mu, R}:=\sup _{|y| \geq R}\left(e^{\mu|y|}|f(y)|\right)$.
\end{definition}

\begin{lemma}[Exponentially Weighted Sup-Norm Estimate]\label{lem:Exponentially Weighted Sup-Norm Estimate}
Let $\Pi$ be a Euclidean space. Let $f \in C^0(\Pi, \mathbb{R}_{\ge 0})$ satisfy the a priori limit $\lim_{|y| \to \infty} f(y) = 0$.

Let $\mu > 0$. Let $E \in \mathcal{B}_\infty(e^{-\mu|y|})$, strictly implying there exists a constant $C_E > 0$ such that $e^{\mu|y|}E(y) \le C_E$ for all sufficiently large $|y|$.

Let $K: \Pi \times \Pi \to \mathbb{R}_{\ge 0}$ be a continuous integration kernel satisfying the uniform weighted bound: 
\begin{equation}
\sup_{y_0} \int_{B_{|y_0|/2}(y_0)} K(y_0, \zeta)e^{\mu|y_0 - \zeta|} d\zeta \le C_K < \infty.
\end{equation}
Assume that for all sufficiently large $y_0 \in \Pi$, $f$ satisfies the integral inequality:
\begin{equation}
f(y_0) \le \int_{B_{|y_0|/2}(y_0)} K(y_0, \zeta) f(\zeta)^3 d\zeta + E(y_0).
\end{equation}
Then, there exists a critical radius $R_{crit} > 0$ such that the weighted norm, defined as
\begin{equation}
\|f\|_{\mu, R} := \sup_{|y| \ge R} (e^{\mu|y|}f(y)),
\end{equation}
 satisfies
\begin{equation}
\|f\|_{\mu, R} \le 2C_E \quad \text{for all} \quad R \ge R_{crit}.
\end{equation}
\end{lemma}

\begin{proof}
Let $y_0 \in \Pi$ be sufficiently large and define $R = |y_0|$. We multiply the given integral inequality by the strictly positive exponential weight $e^{\mu R}$:
\begin{equation}
e^{\mu|y_0|}f(y_0) \le \int_{B_{R/2}(y_0)} K(y_0, \zeta)e^{\mu|y_0|} f(\zeta)^3 d\zeta + e^{\mu|y_0|}E(y_0).
\end{equation}
By the metric triangle inequality on the Euclidean space $\Pi$, we have $|y_0| = |(y_0 - \zeta) + \zeta| \le |y_0 - \zeta| + |\zeta|$. Because the natural exponential function is strictly monotonically increasing over $\mathbb{R}$, we deduce:
\begin{equation}
e^{\mu|y_0|} \le e^{\mu|y_0 - \zeta|}e^{\mu|\zeta|}.
\end{equation}
Substituting this into the integrand and algebraically factoring the cubic function yields:
\begin{equation}
e^{\mu|y_0|}f(y_0) \le \int_{B_{R/2}(y_0)} K(y_0, \zeta)e^{\mu|y_0 - \zeta|} (e^{\mu|\zeta|}f(\zeta)) f(\zeta)^2 d\zeta + e^{\mu|y_0|}E(y_0).
\end{equation}
For any integration variable $\zeta \in B_{R/2}(y_0)$, the reverse triangle inequality imposes the strict constraint $|\zeta| \ge |y_0| - |y_0 - \zeta| \ge R - R/2 = R/2$. Consequently, we can extract the strict algebraic suprema over the exterior domain $\{ \zeta : |\zeta| \ge R/2 \}$ outside of the integral operator:
\begin{equation}
e^{\mu|\zeta|}f(\zeta) \le \sup_{|x| \ge R/2} (e^{\mu|x|}f(x)) = \|f\|_{\mu, R/2},
\end{equation}
\begin{equation}
f(\zeta)^2 \le \sup_{|x| \ge R/2} f(x)^2 = g(R/2)^2,
\end{equation}
where $g(r) := \sup_{|x| \ge r} f(x)$.

Factoring these scalar suprema out of the integral gives:
\begin{equation}
e^{\mu|y_0|}f(y_0) \le g(R/2)^2 \|f\|_{\mu, R/2} \int_{B_{R/2}(y_0)} K(y_0, \zeta)e^{\mu|y_0 - \zeta|} d\zeta + e^{\mu|y_0|}E(y_0).
\end{equation}
By assumption and defining the pointwise weighted error $\tilde{E}(y_0) = e^{\mu|y_0|}E(y_0)$, we obtain:
\begin{equation}
e^{\mu|y_0|}f(y_0) \le C_K g(R/2)^2 \|f\|_{\mu, R/2} + \tilde{E}(y_0).
\end{equation}
Taking the strict supremum over all $y_0$ satisfying $|y_0| \ge R$ yields the global functional inequality:
\begin{equation}
\|f\|_{\mu, R} \le C_K g(R/2)^2 \|f\|_{\mu, R/2} + \|\tilde{E}\|_{\infty, R}.
\end{equation}
By assumption, ($\lim_{|y| \to \infty} f(y) = 0$), the bounding function $g(R/2) \to 0$ as $R \to \infty$. Thus, there exists a radius $R_{crit} > 0$ such that for all $R \ge R_{crit}$, the coefficient strictly satisfies $C_K g(R/2)^2 \le 1/2$.

Furthermore, by assumption, the weighted error norm $\|\tilde{E}\|_{\infty, R}$ is globally bounded by the constant $C_E$ for all $R \ge R_{crit}$. Substituting these bounds, the functional iteratively satisfies:
\begin{equation}
\|f\|_{\mu, R} \le \frac{1}{2} \|f\|_{\mu, R/2} + C_E.
\end{equation}

By defining the remainder function as a strict constant $B(R) = C_E$, the geometric scale ratio becomes exactly $M = B(R/2)/B(R) = 1$. The contraction constant is $\theta = 1/2$. Because $\theta M = 1/2 < 1$ , Lemma~\ref{lem:Iterative Absorption}  absorbs the weighted functional:
\begin{equation}
\|f\|_{\mu, R} \le \frac{1}{1 - 1/2} C_E = 2C_E.
\end{equation}
This completes the proof.
\end{proof}

\begin{remark}
This lemma formalizes a non-linear absorption principle. By leveraging the strictly localized nature of the continuous integration kernel $K$, the non-linear volume contribution is iteratively absorbed, demonstrating that the function $f$ canonically inherits the exponential decay of the a priori error term $E$.
\end{remark}

\section{Asymptotic Analysis and Linearization}

Having established the structural foundations, we proceed to the core analytic exposition of the translator equation. This chapter is devoted to the systematic reduction of the quasilinear mean curvature condition into a tractable, asymptotically linear form. By extracting the principal linear differential operator and rigorously bounding the residual nonlinearities, we establish the global pointwise decay estimates necessary to enforce geometric rigidity.

\subsection{Quasilinear Elliptic Structure and Geometric Reductions}

We first translate the purely geometric translator condition into a strict analytic formalism, obtaining a quasilinear elliptic partial differential equation governing the scalar graph parametrization $u$. We establish the fundamental invariance of this differential operator under vertical rigid translations. This geometric reduction is analytically significant, as it permits the canonical normalization of the asymptotic end, strictly forcing the parameterization to vanish at infinity ($\lim_{|y|\to\infty} |\tilde{u}(y)| = 0$), thereby simplifying all subsequent global estimates.

\begin{lemma}[Quasilinear Translator Equation]\label{lem:Quasilinear Translator Equation}
Let $\Sigma$ be locally represented as the graph of $u: \Omega \subset \Pi \to \mathbb{R}$. The geometric condition $\mathbf{H} = e_{n+1}^\perp$ is equivalent to the quasilinear elliptic equation:
\begin{equation}
\mathcal{L}[u] := \sum_{i,j=1}^n \left( \delta_{ij} - \frac{D_i u D_j u}{1+|Du|^2} \right) D_{ij} u + D_n u = 0
\end{equation}
where $D$ denotes the gradient on $\Pi \cong \mathbb{R}^n$.
\end{lemma}

\begin{proof}
We fix the ambient Euclidean space $\mathbb{R}^{n+1}$ equipped with the standard Cartesian coordinates $X = (x_1, \dots, x_{n+1})$ and the orthonormal basis $\{e_1, \dots, e_{n+1}\}$.  We designate $e_{n+1}$ as the translation direction. Let $\Pi$ denote the vertical hyperplane normalized to be the coordinate plane $\{x_1 = 0\}$. We adopt coordinates $y = (y_1, \dots, y_n)$ on $\Pi$, identifying the $y_n$-axis with the $x_{n+1}$-direction.

The surface $\Sigma$ is the image of the embedding $F: \Omega \to \mathbb{R}^{n+1}$ defined by $F(y) = (u(y), y_1, \dots, y_n)$. To investigate the intrinsic geometry, we compute the tangent basis vectors $E_i$ on $\Sigma$. Differentiating the position vector with respect to the coordinates $y_i$ yields $E_i = D_i u \, e_1 + \hat{e}_i$, where $\hat{e}_i$ denotes the basis vector in the $y_i$-direction within the hyperplane $\Pi$. The components of the induced metric tensor $g_{ij}$ are given by the inner product
\begin{equation}
g_{ij} = \langle E_i, E_j \rangle = \delta_{ij} + D_i u D_j u
\end{equation}
where $D$ denotes the gradient on $\Pi$. The inverse metric tensor, denoted by $g^{ij}$, is obtained by applying the Sherman-Morrison formula to the rank-one perturbation of the identity matrix. A straightforward computation yields
\begin{equation}
g^{ij} = \delta_{ij} - \frac{D_i u D_j u}{1 + |Du|^2}
\end{equation}
where $|Du|^2 = \sum (D_k u)^2$.

We proceed to determine the unit normal vector field $\nu$ along $\Sigma$. The vector $\nu$ must be orthogonal to the tangent space spanned by $\{E_i\}$. Consider the vector $N = e_1 - \sum_{k=1}^n (D_k u) \hat{e}_k$. Normalizing $N$, we obtain the unit normal
\begin{equation}
\nu = \frac{e_1 - \sum_{k=1}^n (D_k u) \hat{e}_k}{\sqrt{1 + |Du|^2}}
\end{equation}
We calculate the coefficients of the second fundamental form, $h_{ij}$. The second derivatives of the position vector are strictly in the $e_1$ direction, specifically $D_{ij} F = D_{ij} u \, e_1$. Consequently, the second fundamental form is
\begin{equation}
h_{ij} = \langle D_{ij} F, \nu \rangle = \frac{D_{ij} u}{\sqrt{1 + |Du|^2}}
\end{equation}
The mean curvature scalar $H$ is the trace of the second fundamental form with respect to the induced metric. Thus, we have
\begin{equation}
H = g^{ij} h_{ij} = \sum_{i,j=1}^n \left( \delta_{ij} - \frac{D_i u D_j u}{1 + |Du|^2} \right) \frac{D_{ij} u}{\sqrt{1 + |Du|^2}}
\end{equation}

The geometric definition of a translator states that the mean curvature vector satisfies $\mathbf{H} = e_{n+1}^\perp$. Expressing this in terms of the scalar mean curvature relative to $\nu$, we have $H = \langle e_{n+1}, \nu \rangle$. We must evaluate the inner product on the right-hand side. Recalling that the ambient translation direction $e_{n+1}$ corresponds to the coordinate axis $y_n$ on $\Pi$, its projection onto the normal $\nu$ is
\begin{equation}
\langle e_{n+1}, \nu \rangle = - \frac{D_n u}{\sqrt{1 + |Du|^2}}
\end{equation}
Equating the two expressions for $H$, we obtain
\begin{equation}
\sum_{i,j=1}^n \left( \delta_{ij} - \frac{D_i u D_j u}{1 + |Du|^2} \right) \frac{D_{ij} u}{\sqrt{1 + |Du|^2}} = - \frac{D_n u}{\sqrt{1 + |Du|^2}}
\end{equation}
Since the factor $(1 + |Du|^2)^{-1/2}$ is strictly positive, we may multiply the entire equation by this term. Rearranging terms to one side, we arrive at the requisite quasilinear elliptic equation
\begin{equation}
\sum_{i,j=1}^n \left( \delta_{ij} - \frac{D_i u D_j u}{1 + |Du|^2} \right) D_{ij} u + D_n u = 0
\end{equation}
This concludes the proof.
\end{proof}

\begin{remark}
The passage from the extrinsic geometric relation of the immersed manifold to an intrinsic partial differential equation is canonical. By expanding the mean curvature as the trace of the second fundamental form with respect to the induced metric, the geometric constraint is strictly reduced to a quasilinear structure, thereby permitting the deployment of standard elliptic regularity machinery.
\end{remark}

\begin{proposition}[Vertical Translation Invariance]\label{prop:Vertical Translation Invariance}
Let $\Omega \subset \mathbb{R}^n$ be an open domain. Let $\mathcal{L}$ be the quasilinear operator defined by
\begin{equation}
\mathcal{L}[u] := \sum_{i,j=1}^{n} \left(\delta_{ij} - \frac{D_i u D_j u}{1 + |D u|^2}\right) D_{ij} u + D_n u.
\end{equation}
For any constant $c \in \mathbb{R}$ and any function $u \in C^2(\Omega, \mathbb{R})$, let $v$ be the translated function defined pointwise by $v(y) = u(y) - c$. Then $\mathcal{L}[v](y) = \mathcal{L}[u](y)$ for all $y \in \Omega$.
\end{proposition}

\begin{proof}
Let $y \in \Omega$. We evaluate the function $v(y) = u(y) - c$.

First, evaluate the first derivative. By the linearity of the gradient operator,we have
\begin{equation}
D v(y) = D(u(y) - c) = D u(y) - D(c).
\end{equation}
The gradient of the constant field $c$ is the zero vector. Therefore,
\begin{equation}
D v(y) = D u(y) - \mathbf{0} = D u(y),
\end{equation}
which implies that for every component $i \in \{1, \dots, n\}$, $D_i v(y) = D_i u(y)$, and thus $|D v(y)|^2 = |D u(y)|^2$.

Next, evaluate the second derivative. By the linearity of the Hessian operator, we get
\begin{equation}
D^2 v(y) = D^2(u(y) - c) = D^2 u(y) - D^2(c).
\end{equation}
Then $D^2(c) = \mathbf{0}$. Therefore,
\begin{equation}
D^2 v(y) = D^2 u(y) - \mathbf{0} = D^2 u(y),
\end{equation}
which implies that for every matrix component $i,j \in \{1, \dots, n\}$, $D_{ij} v(y) = D_{ij} u(y)$.

Now, substitute into the operator expression. We have
\begin{equation}
\mathcal{L}[v](y) = \sum_{i,j=1}^{n} \left(\delta_{ij} - \frac{D_i v(y) D_j v(y)}{1 + |D v(y)|^2}\right) D_{ij} v(y) + D_n v(y).
\end{equation}
Substituting the established equalities $D v(y) = D u(y)$ and $D^2 v(y) = D^2 u(y)$ into this expression, we get
\begin{equation}
\mathcal{L}[v](y) = \sum_{i,j=1}^{n} \left(\delta_{ij} - \frac{D_i u(y) D_j u(y)}{1 + |D u(y)|^2}\right) D_{ij} u(y) + D_n u(y).
\end{equation}
Thus, we find that $\mathcal{L}[v](y)$ is strictly equal to $\mathcal{L}[u](y)$. This completes the proof.
\end{proof}

\begin{corollary}[Asymptotic Normalization]\label{cor:Asymptotic Normalization}
Let $E$ be an end $C^1$-asymptotic to a vertical hyperplane $\Pi$, represented by $graph(u)$ with $\lim_{|y| \to \infty} ||Du(y)|| = 0$. There exists a rigid translation of $\Sigma$ along the $e_1$ normal vector such that the resulting parameterization $\tilde{u}$  satisfies both $\lim_{|y| \to \infty} ||D\tilde{u}(y)|| = 0$ and $\lim_{|y| \to \infty} |\tilde{u}(y)| = 0$.
\end{corollary}

\begin{proof}
Let $u(y)$ be the parameterization of the end $E$ strictly satisfying $\lim_{|y| \to \infty} ||Du(y)|| = 0$. Then there exists a finite real constant $c \in \mathbb{R}$ such that $\lim_{|y| \to \infty} u(y) = c$.

We apply a geometric rigid translation to $\Sigma$ by the vector $-c e_1$. This strictly induces the new algebraic parameterization $\tilde{u}(y) = u(y) - c$.

By the linearity of limits, we evaluate the limit of the translated parameterization as $|y| \to \infty$:
\begin{equation}
\lim_{|y| \to \infty} \tilde{u}(y) = \lim_{|y| \to \infty} (u(y) - c) = \left( \lim_{|y| \to \infty} u(y) \right) - c.
\end{equation}
Substituting,
\begin{equation}
\lim_{|y| \to \infty} \tilde{u}(y) = c - c = 0.
\end{equation}
Taking the absolute value  yields
\begin{equation}
\lim_{|y| \to \infty} |\tilde{u}(y)| = 0.
\end{equation}

By the linearity of the gradient operator $D$, we compute the pointwise gradient of $\tilde{u}$:
\begin{equation}
D\tilde{u}(y) = D(u(y) - c) = Du(y) - D(c).
\end{equation}
Because the derivative of a constant scalar field is the zero vector ($\mathbf{0}$):
\begin{equation}
D\tilde{u}(y) = Du(y) - \mathbf{0} = Du(y).
\end{equation}
Therefore, taking the limit of the vector norm yields:
\begin{equation}
\lim_{|y| \to \infty} ||D\tilde{u}(y)|| = \lim_{|y| \to \infty} ||Du(y)||.
\end{equation}
By the initial hypothesis, this strictly evaluates to $0$.

The rigid translation $v(y) = u(y) - c$  results in a parameterization $\tilde{u}$ that satisfies both $\lim_{|y| \to \infty} ||D\tilde{u}(y)|| = 0$ and $\lim_{|y| \to \infty} |\tilde{u}(y)| = 0$. 
\end{proof}

\begin{remark}
Modulo the natural action of the affine group of translations along the $e_1$ normal vector in the ambient Euclidean space, one may always canonically normalize the asymptotic datum to zero. This rigid translation strictly preserves the triviality of the limiting gradient.
\end{remark}

\subsection{Asymptotic Analysis and the Drift Laplacian}

Under the hypothesis of asymptotic gradient triviality, the quasilinear operator exhibits a well-defined limiting structure. We prove that the operator is uniformly elliptic on the exterior domain and structurally converges to the constant-coefficient linear drift Laplacian, denoted $L_0 := \Delta + \partial_{y_n}$. We subsequently derive strict decay estimates for the fundamental solution of this specific drift Laplacian and construct the formal Green representation formula for the pointwise gradient. This effectively linearizes the asymptotic geometry, isolating the translation drift as the dominant structural feature.

\begin{definition}(Asymptotic Translating Graph)
Let $E$ be a real Euclidean space of dimension $n$. Let $\Omega \subset E$ be an exterior domain. A function $u \in C^2(\Omega)$ is an asymptotic translating graph if it satisfies the quasilinear equation:

\begin{equation} 
\sum_{i, j=1}^n\left(\delta_{i j}-\frac{D_i u D_j u}{1+|D u|^2}\right) D_{i j} u+D_n u=0
\end{equation}

and strictly satisfies the asymptotic gradient triviality condition $\lim _{|y| \rightarrow \infty}|D u(y)|=0$.
\end{definition}

\begin{definition}(Drift Laplacian)
Let $E$ be a real Euclidean space of dimension $n$. The standard drift Laplacian is the linear operator $L_0: C^2(E) \rightarrow C^0(E)$ defined by $L_0[v]=\Delta v+\partial_{y_n} v$.
\end{definition}

\begin{proposition}[Structure of the Linearized Operator]\label{prop:Structure of the Linearized Operator}
Let $u$ be a solution to $\mathcal{L}[u] = 0$ on an exterior domain $\Pi \setminus B_R$. If $\lim_{|y| \to \infty} |Du| = 0$, then $\mathcal{L}$ is uniformly elliptic and asymptotically equivalent to the linear Drift Laplacian:
\begin{equation}
L_0 := \Delta + \partial_{y_n}
\end{equation}
\end{proposition}

\begin{proof}
Let $u$ be an Asymptotic Translating Graph. Then its associated operator $\mathcal{L}$ is asymptotically equivalent to the Linearized Drift Operator $L_0$. Let us write the quasilinear operator $\mathcal{L}$ in non-divergence form. The equation is given by
\begin{equation}
\mathcal{L}[u] = \sum_{i,j=1}^n a_{ij}(Du) D_{ij} u + b(Du) \cdot Du
\end{equation}
where the coefficient matrix $A = (a_{ij})$ and the vector field term are defined, respectively, by the structural components of the mean curvature operator and the translation projection. Explicitly, we have
\begin{equation}
a_{ij}(p) = \delta_{ij} - \frac{p_i p_j}{1+|p|^2} \quad \text{and} \quad \text{Drift term } D_n u
\end{equation}
where $p = (p_1, \dots, p_n)$ represents the gradient $Du$. Note that the drift term $D_n u$ corresponds to a linear first-order operator $\partial_{y_n}$ acting on $u$.

We proceed by verifying the ellipticity of the coefficient matrix $A$. For an arbitrary non-zero vector $\xi \in \mathbb{R}^n$, we evaluate the quadratic form associated with $a_{ij}$:
\begin{equation}
Q_A(\xi) := \sum_{i,j=1}^n a_{ij}(p) \xi_i \xi_j = \sum_{i,j=1}^n \left( \delta_{ij} - \frac{p_i p_j}{1+|p|^2} \right) \xi_i \xi_j
\end{equation}
Expanding the terms linearly, we obtain
\begin{equation}
Q_A(\xi) = |\xi|^2 - \frac{(\sum_{i=1}^n p_i \xi_i)^2}{1+|p|^2} = |\xi|^2 - \frac{\langle p, \xi \rangle^2}{1+|p|^2}
\end{equation}
By the Cauchy-Schwarz inequality, we have $\langle p, \xi \rangle^2 \le |p|^2 |\xi|^2$. This allows us to bound the quadratic form from below:
\begin{equation}
Q_A(\xi) \ge |\xi|^2 - \frac{|p|^2 |\xi|^2}{1+|p|^2} = |\xi|^2 \left( 1 - \frac{|p|^2}{1+|p|^2} \right) = \frac{|\xi|^2}{1+|p|^2}
\end{equation}
Conversely, since the subtractive term is non-negative, we immediately have the upper bound $Q_A(\xi) \le |\xi|^2$. Thus, the eigenvalues $\lambda(A)$ of the matrix $A$ lie strictly within the interval
\begin{equation}
\frac{1}{1+|p|^2} \le \lambda(A) \le 1
\end{equation}

The assumption says that the gradient decays at infinity, specifically $\lim_{|y| \to \infty} |Du(y)| = 0$. Consequently, for any $\epsilon > 0$, there exists a radius $R_\epsilon \ge R$ such that for all $y \in \Pi \setminus B_{R_\epsilon}$, we have $|Du(y)| < \epsilon$. Choosing $\epsilon$ sufficiently small, we observe that the ellipticity constants are uniformly bounded away from zero and infinity on the exterior domain. Specifically, we have
\begin{equation}
\frac{1}{2} \le \frac{1}{1+|Du|^2} \le \lambda(A) \le 1
\end{equation}
This confirms that $\mathcal{L}$ acts as a uniformly elliptic operator on the domain $\Pi \setminus B_{R_\epsilon}$.

Finally, we determine the asymptotic limit of the operator structure as $|y| \to \infty$. We examine the behavior of the coefficients $a_{ij}(Du)$ in the limit of vanishing gradient:
\begin{equation}
\lim_{|y| \to \infty} a_{ij}(Du(y)) = \lim_{|p| \to 0} \left( \delta_{ij} - \frac{p_i p_j}{1+|p|^2} \right) = \delta_{ij} - 0 = \delta_{ij}
\end{equation}
The linear part of the operator is strictly $D_n$, which is $\partial_{y_n}$. Therefore, in the asymptotic limit, the operator $\mathcal{L}$ acting on $u$ converges structurally to the constant-coefficient operator composed of the limiting principal part and the drift term:
\begin{equation}
\mathcal{L}_{asymp} = \sum_{i,j=1}^n \delta_{ij} D_{ij} + \partial_{y_n} = \Delta + \partial_{y_n}
\end{equation}
We define this limit as $L_0 := \Delta + \partial_{y_n}$. Thus, $\mathcal{L}$ is asymptotically equivalent to the linear Drift Laplacian $L_0$.
\end{proof}

\begin{proposition}[Fundamental Solution Decay Estimate]\label{prop:Fundamental Solution Decay Estimate}
Let $E = \mathbb{R}^n$ with $n \ge 3$, and let $L_0 : C^2(E) \to C^0(E)$ be the drift Laplacian defined by $L_0 = \Delta + \partial_{y_n}$. Let $K \in C^2(E \setminus {0})$ be the fundamental solution satisfying $L_0 K = \delta_0$ in the sense of distributions. Then there exists a constant $C>0$ such that for all $|y| \ge 1$ one has
\begin{equation} 
|K(y)| \le C |y|^{-(n-1)/2} e^{-\frac12(|y|+y_n)} .
\end{equation}
\end{proposition}

\begin{proof}
The fundamental solution of the heat equation $\partial_t u = \Delta u$ on $\mathbb{R}^n$ is the Gaussian kernel $p(t,y) = (4\pi t)^{-n/2}\exp(-|y|^2/(4t))$. By Galilean invariance, the fundamental solution of the drift equation $\partial_t u = (\Delta+\partial_{y_n})u$ is obtained by a translation in the $-e_n$ direction, and is given by
\begin{equation} 
p_D(t,y) = (4\pi t)^{-n/2}\exp!\left(-\frac{|y+te_n|^2}{4t}\right).
\end{equation}

The stationary fundamental solution of $L_0$ is obtained by integrating the negative parabolic kernel over $t>0$, so
\begin{equation} 
K(y) = -\int_0^\infty (4\pi t)^{-n/2}\exp!\left(-\frac{|y+te_n|^2}{4t}\right),dt .
\end{equation}

Expanding $|y+te_n|^2 = |y|^2 + 2ty_n + t^2$ gives
\begin{equation} 
-\frac{|y+te_n|^2}{4t} = -\frac{|y|^2}{4t} - \frac{y_n}{2} - \frac{t}{4}.
\end{equation}
Factoring the constant exponential term yields
\begin{equation} 
|K(y)| = e^{-y_n/2}\int_0^\infty (4\pi t)^{-n/2}\exp!\left(-\frac{|y|^2}{4t}-\frac{t}{4}\right),dt .
\end{equation}

Denote the integral factor by $I(y)$. Introduce the change of variables $t = |y|s$ with $dt = |y|ds$. Substituting gives
\begin{equation} 
I(y) = \frac{|y|^{1-n/2}}{(4\pi)^{n/2}}\int_0^\infty s^{-n/2}\exp!\left(-\frac{|y|}{4}\left(s+\frac1s\right)\right),ds .
\end{equation}

Using the identity $s + \frac1s = 2 + \frac{(s-1)^2}{s}$ we obtain
\begin{equation} 
I(y) = \frac{|y|^{1-n/2}}{(4\pi)^{n/2}} e^{-|y|/2}\int_0^\infty s^{-n/2}\exp!\left(-\frac{|y|}{4}\frac{(s-1)^2}{s}\right),ds .
\end{equation}

Let $J(|y|) = \int_0^\infty s^{-n/2}\exp!\left(-\frac{|y|}{4}\frac{(s-1)^2}{s}\right)ds$. For $|y|\ge1$ we partition the domain $(0,\infty)$ into $(0,1/2)$, $[1/2,2]$, and $(2,\infty)$. On $[1/2,2]$ the factor $s^{-n/2}$ is bounded and $1/s \ge 1/2$, so the exponential term produces a Gaussian bound and the integral over this interval is bounded by $C_1 |y|^{-1/2}$. On $(2,\infty)$ one has $\frac{(s-1)^2}{s} > s/4$, so the exponent is bounded above by $-|y|s/16$, which yields $\int_2^\infty \cdots ds \le C_2 |y|^{-1/2}$. On $(0,1/2)$ we apply the substitution $v = 1/s$; the symmetry of the integrand under this inversion gives $\int_0^{1/2} \cdots ds \le C_3 |y|^{-1/2}$. Combining the three bounds yields $J(|y|) \le C_0 |y|^{-1/2}$ for all $|y| \ge 1$.

Substituting this estimate into the expression for $I(y)$ gives $I(y) \le C_4 |y|^{-(n-1)/2} e^{-|y|/2}$. Multiplying by the previously extracted factor $e^{-y_n/2}$ yields
\begin{equation} 
|K(y)| \le C_4 |y|^{-(n-1)/2} e^{-\frac12(|y|+y_n)} .
\end{equation}
\end{proof}

\begin{proposition}[Gradient Green Representation]\label{prop:Gradient Green Representation}
Let $\Omega \subset \mathbb{R}^n$ with $n \ge 3$ be an open set. Let $u \in C^2(\Omega)$ satisfy $L_0 u = f$ in $\Omega$, where $L_0 = \Delta + \partial_{y_n}$ and $f \in C^0(\Omega)$. Let $y_0 \in \Omega$ be such that the closed ball $\overline{B_\rho(y_0)} \subset \Omega$. Let $K \in C^2(\mathbb{R}^n \setminus {0})$ be the fundamental solution satisfying $L_0 K = \delta_0$. Then the pointwise gradient satisfies
\begin{equation} 
|Du(y_0)| \le \int_{B_\rho(y_0)} |D_y K(y_0-\zeta)|, |f(\zeta)|, d\zeta + \int_{\partial B_\rho(y_0)} \mathcal{R}(u,K,\zeta), dS(\zeta),
\end{equation}
where the boundary remainder density is
\begin{equation} 
\mathcal{R}(u,K,\zeta) = |u(\zeta)|, |D_y \partial_\nu K(y_0-\zeta)| + |\partial_\nu u(\zeta)|, |D_y K(y_0-\zeta)| + |u(\zeta)|, |\nu_n(\zeta)|, |D_y K(y_0-\zeta)|.
\end{equation}
\end{proposition}

\begin{proof}
Let $L_0^* v = \Delta v - \partial_{\zeta_n} v$ denote the formal adjoint of $L_0$. Consider the vector field $W = v\nabla u - u\nabla v + uv e_n$. A direct computation gives $\nabla \cdot W = v L_0 u - u L_0^* v$. Let $v(\zeta) = K(y_0-\zeta)$. Applying the divergence theorem on the punctured domain $V_\epsilon = B_\rho(y_0) \setminus \overline{B_\epsilon(y_0)}$ and passing to the limit as $\epsilon \to 0$ yields the Green representation formula
\begin{equation} 
u(y_0) = \int_{B_\rho(y_0)} K(y_0-\zeta) f(\zeta), d\zeta - \int_{\partial B_\rho(y_0)} B(u(\zeta),K(y_0-\zeta)), dS(\zeta),
\end{equation}
where the boundary operator is $B(u,K) = K \partial_\nu u - u \partial_\nu K + uK \nu_n$.

We differentiate this identity with respect to the parameter $y_0$. Consider first the volume integral. Let the difference quotient be $D^h(\zeta) = \frac{1}{h}\big(K(y_0+he_i-\zeta)-K(y_0-\zeta)\big)f(\zeta)$. By the mean value theorem one has $|D^h(\zeta)| \le C |y_0+\theta he_i-\zeta|^{1-n}\sup |f|$ for some $\theta \in (0,1)$. Because $|x|^{1-n}$ belongs to $L^1(B_\rho)$ when $n \ge 3$, the Lebesgue dominated convergence theorem implies
\begin{equation} 
D_{y_0}\int_{B_\rho(y_0)} K(y_0-\zeta) f(\zeta), d\zeta = \int_{B_\rho(y_0)} D_{y_0} K(y_0-\zeta) f(\zeta), d\zeta .
\end{equation}

Next we differentiate the boundary integral. On $\partial B_\rho(y_0)$ the distance between $y_0$ and $\zeta$ equals $\rho$, so the kernel and its derivatives are smooth. The quantities $u(\zeta)$, $\partial_\nu u(\zeta)$ and $\nu_n(\zeta)$ depend only on $\zeta$ and are independent of $y_0$. Therefore the derivative acts only on the kernel. Expanding the derivative gives
\begin{equation} 
D_{y_0}[B(u,K(y_0-\zeta))] = D_{y_0}K(y_0-\zeta),\partial_\nu u(\zeta) - u(\zeta), D_{y_0}[\partial_\nu K(y_0-\zeta)] + u(\zeta), D_{y_0}K(y_0-\zeta),\nu_n(\zeta).
\end{equation}

By translation invariance one has $D_{y_0}K(y_0-\zeta) = -D_\zeta K(y_0-\zeta)$, which we denote by $-D_y K(y_0-\zeta)$. Taking absolute values and applying the triangle inequality to the differentiated representation formula yields
\begin{equation} 
|Du(y_0)| \le \int_{B_\rho(y_0)} |D_y K(y_0-\zeta)|, |f(\zeta)|, d\zeta
* \int_{\partial B_\rho(y_0)} \left| (D_y K)\partial_\nu u - u(D_y\partial_\nu K) + u(D_y K)\nu_n \right| dS(\zeta).
\end{equation}

Distributing the absolute value across the boundary terms gives
\begin{equation} 
\mathcal{R}(u,K,\zeta) = |\partial_\nu u(\zeta)|, |D_y K(y_0-\zeta)| + |u(\zeta)|, |D_y \partial_\nu K(y_0-\zeta)| + |u(\zeta)|, |\nu_n(\zeta)|, |D_y K(y_0-\zeta)|.
\end{equation}
Substituting this expression completes the estimate.
\end{proof}

\begin{remark}
The pointwise control of the gradient vector field is thus canonically reconstructed through the Green's formalism. This mapping strategically isolates the non-linearities into a volume convolution, while strictly encapsulating the far-field geometric data within a boundary remainder density over the sphere.
\end{remark}

\subsection{Nonlinear Remainder and Global Decay Theorems}

The final analytic requirement is the strict domination of the nonlinear remainder operator, defined exactly as the formal algebraic difference between the quasilinear operator and the linear drift Laplacian ($\mathcal{N}[u] := \mathcal{L}[u] - L_0[u]$). We establish explicit pointwise differential inequalities bounding this remainder by the Hessian and squared gradient. By constructing an $\epsilon$-tilted formal barrier function and invoking the comparison principle for exterior domains, we derive an a priori $C^0$ global decay estimate. This machinery culminates in the proof of bifurcated, pointwise decay rates for the gradient vector field, structurally segregating its behavior across the aforementioned wake and exponential sets.

\begin{definition}(symmetric Difference Operator)
Let $\Omega \subset \mathbb{R}^n$ be an open domain. For any $u, V \in C^2(\Omega, \mathbb{R})$, define the difference function $v=u-V$. We define the canonical phase paths for the gradient and the Hessian as:

\begin{equation} 
\begin{aligned}
\gamma_1(t, y) & =t D u(y)+(1-t) D V(y) \\
\gamma_2(t, y) & =t D^2 u(y)+(1-t) D^2 V(y)
\end{aligned}
\end{equation}

We define the coefficient matrix $\tilde{a}^{i j}(y)$ and the drift vector field $b^k(y)$ via the continuous integrals:

\begin{equation} 
\begin{array}{c}
\tilde{a}^{i j}(y)=\int_0^1 a^{i j}\left(\gamma_1(t, y)\right) d t \\
b^k(y)=\int_0^1 \frac{\partial a^{i j}}{\partial p_k}\left(\gamma_1(t, y)\right)\left[\gamma_2(t, y)\right]_{i j} d t
\end{array}
\end{equation}

where $a^{i j}(p)=\delta_{i j}-\frac{p_i p_j}{1+|p|^2}$.
\end{definition}

\begin{lemma}[Symmetric Difference Linearization]\label{lem:Symmetric Difference Linearization}
If $u, V \in C^2(\Omega, \mathbb{R})$, then the difference of the quasilinear operators exactly satisfies the linear partial differential equation:
\begin{equation}
\mathcal{L}[u] - \mathcal{L}[V] = \tilde{a}^{ij}(y) D_{ij} v + b^k(y) D_k v + D_n v
\end{equation}
\end{lemma}

\begin{proof}
Let $y \in \Omega$ be fixed. The quasilinear operator is defined structurally as $\mathcal{L}[w] = a^{ij}(Dw) D_{ij} w + D_n w$.

We evaluate the strict algebraic difference evaluated at $y$:
\begin{equation}
\mathcal{L}[u] - \mathcal{L}[V] = [a^{ij}(Du) D_{ij} u - a^{ij}(DV) D_{ij} V] + [D_n u - D_n V]
\end{equation}
By the linearity of the differential operator $D_n$, the drift component trivially yields:
\begin{equation}
D_n u - D_n V = D_n (u - V) = D_n v
\end{equation}
We now isolate the principal part. Define the scalar functional $\Phi: [0,1] \to \mathbb{R}$ representing the tensor contraction along the integral paths established in the definition,
\begin{equation}
\Phi(t) = a^{ij}(\gamma_1(t)) [\gamma_2(t)]_{ij}
\end{equation}
Evaluating the boundaries of the interval $[0,1]$ yields:
\begin{equation}
\Phi(1) = a^{ij}(Du) D_{ij} u
\end{equation}
\begin{equation}
\Phi(0) = a^{ij}(DV) D_{ij} V
\end{equation}
Therefore, the principal difference is exactly $\Phi(1) - \Phi(0)$.

Because $u, V \in C^2$, the paths $\gamma_1$ and $\gamma_2$ are strictly continuously differentiable with respect to $t$. By the Fundamental Theorem of Calculus on $\mathbb{R}$:
\begin{equation}
\Phi(1) - \Phi(0) = \int_0^1 \frac{d\Phi}{dt} dt
\end{equation}
We compute the total derivative $\frac{d\Phi}{dt}$ using the multivariable chain rule and the product rule for tensor contractions:
\begin{equation}
\frac{d\Phi}{dt} = \frac{d}{dt}[a^{ij}(\gamma_1(t))] [\gamma_2(t)]_{ij} + a^{ij}(\gamma_1(t)) \frac{d}{dt} [ [\gamma_2(t)]_{ij}]
\end{equation}
Applying the chain rule to the matrix component:
\begin{equation}
\frac{d}{dt}[a^{ij}(\gamma_1(t))] = \sum_{k=1}^n \frac{\partial a^{ij}}{\partial p_k}(\gamma_1(t)) \frac{d}{dt} [\gamma_1(t)]_k
\end{equation}
We compute the phase velocities (the derivatives with respect to $t$):
\begin{equation}
\frac{d}{dt} [\gamma_1(t)]_k = D_k u - D_k V = D_k v
\end{equation}
\begin{equation}
\frac{d}{dt} [ [\gamma_2(t)]_{ij}] = D_{ij} u - D_{ij} V = D_{ij} v
\end{equation}
Substituting the phase velocities back into the total derivative yields:
\begin{equation}
\frac{d\Phi}{dt} = \left[ \sum_{k=1}^n \frac{\partial a^{ij}}{\partial p_k}(\gamma_1(t)) [\gamma_2(t)]_{ij} \right] D_k v + \left[ a^{ij}(\gamma_1(t)) \right] D_{ij} v
\end{equation}
Integrating this expression from $t=0$ to $t=1$, and factoring out the spatial derivatives $D_k v$ and $D_{ij} v$ (which are strictly independent of the integration parameter $t$), we obtain:
\begin{equation}
\Phi(1) - \Phi(0) = \left( \int_0^1 a^{ij}(\gamma_1(t)) dt \right) D_{ij} v + \left( \int_0^1 \sum_{k=1}^n \frac{\partial a^{ij}}{\partial p_k}(\gamma_1(t)) [\gamma_2(t)]_{ij} dt \right) D_k v
\end{equation}
By definition, the integrals inside the parentheses map exactly to the canonical coefficients $\tilde{a}^{ij}(y)$ and $b^k(y)$. Thus:
\begin{equation}
a^{ij}(Du) D_{ij} u - a^{ij}(DV) D_{ij} V = \tilde{a}^{ij}(y) D_{ij} v + b^k(y) D_k v
\end{equation}
Adding the linear drift term $D_n v$ back to the equation  seals the identity:
\begin{equation}
\mathcal{L}[u] - \mathcal{L}[V] = \tilde{a}^{ij}(y) D_{ij} v + b^k(y) D_k v + D_n v
\end{equation}
This completes the proof. 
\end{proof}

\begin{proposition}[Comparison Principle for Exterior Domains]\label{prop:Comparison Principle for Exterior Domains}
Let $\Pi \cong \mathbb{R}^n$. Let $u, V \in C^2(\Pi \setminus B_{R_0})$ satisfy $\mathcal{L}[u] \ge \mathcal{L}[V]$. Assume $\lim_{|y|\to\infty} (u(y) - V(y)) \le 0$ and $u(y) \le V(y)$ on $\partial B_{R_0}$. Then $u(y) \le V(y)$ globally on $\Pi \setminus B_{R_0}$.
\end{proposition}

\begin{proof}
Define the difference function $v(y) = u(y) - V(y)$. By the Lemma~\ref{lem:Symmetric Difference Linearization}, $v$ exactly satisfies the linear partial differential inequality:
\begin{equation}
\tilde{L}[v] = \tilde{a}^{ij}(y) D_{ij} v + b^k(y) D_k v + D_n v \ge 0
\end{equation}
everywhere on the exterior domain $\Pi \setminus B_{R_0}$. Notice that the operator $\tilde{L}$ contains no zeroth-order term (the coefficient of $v$ is exactly 0).

We must prove that $v(y) \le 0$ for all $y \in \Pi \setminus B_{R_0}$.

Let $\epsilon > 0$ be an arbitrarily chosen real number. By the strict asymptotic hypothesis, $\limsup_{|y|\to\infty} v(y) \le 0$. Therefore, by the formal definition of a limit at infinity, there exists a sufficiently large radial threshold $R_\epsilon > R_0$ such that for all $y$ satisfying $|y| \ge R_\epsilon$:
\begin{equation}
v(y) \le \epsilon
\end{equation}
We construct the open, bounded, connected annulus $\Omega_\epsilon \subset \mathbb{R}^n$:
\begin{equation}
\Omega_\epsilon = \{y \in \Pi : R_0 < |y| < R_\epsilon\}
\end{equation}
The topological boundary of this annulus is the disjoint union of two spheres: $\partial \Omega_\epsilon = \partial B_{R_0} \cup \partial B_{R_\epsilon}$.

We evaluate the boundary conditions of $v$ on $\partial \Omega_\epsilon$:

On the inner boundary $\partial B_{R_0}$, the hypothesis strictly states $u \le V$, which implies $v(y) \le 0$. Since $\epsilon > 0$, we have $v(y) < \epsilon$.

On the outer boundary $\partial B_{R_\epsilon}$, our choice of the radial threshold guarantees $v(y) \le \epsilon$.

Therefore, over the entire compact boundary $\partial \Omega_\epsilon$, we have the strict bound $\sup_{\partial \Omega_\epsilon} v \le \epsilon$.

By Proposition~\ref{prop:Structure of the Linearized Operator}, the linearized operator $\tilde{L}$ is uniformly elliptic with bounded coefficients over the compact closure $\overline{\Omega_\epsilon}$. Because $v \in C^2(\Omega_\epsilon) \cap C^0(\overline{\Omega_\epsilon})$, $\tilde{L}[v] \ge 0$, and $\tilde{L}$ has a zero coefficient for the zeroth-order term, $v$ is a strict subsolution.

We may now  invoke the classical Weak Maximum Principle for bounded domains, which states that a subsolution to such an elliptic operator achieves its maximum on the boundary:
\begin{equation}
\sup_{\Omega_\epsilon} v \le \sup_{\partial \Omega_\epsilon} v \le \epsilon
\end{equation}
Thus, $v(y) \le \epsilon$ for all $y \in \Omega_\epsilon$.

Furthermore, by the definition of $R_\epsilon$, $v(y) \le \epsilon$ for all $|y| \ge R_\epsilon$.

Combining these domains, we deduce that $v(y) \le \epsilon$ globally for all $y \in \Pi \setminus B_{R_0}$.

Because this global inequality holds strictly for any arbitrary $\epsilon > 0$, we may take the limit as $\epsilon \to 0^+$ to conclude:
\begin{equation}
v(y) \le 0 \quad \text{for all } y \in \Pi \setminus B_{R_0}
\end{equation}
Substituting $v(y) = u(y) - V(y)$, we obtain $u(y) \le V(y)$ globally on $\Pi \setminus B_{R_0}$. This completes the proof.
\end{proof}

\begin{lemma}(Explicit Spherical Anisotropic Averaging)\label{lem:Explicit Spherical Anisotropic Averaging}
Let $n \ge 3$. Let $L_0 = \Delta + \partial_{y_n}$ be the drift Laplacian and let $K(x)$ be its fundamental solution. Assume the kernel's normal derivative satisfies the pointwise bound
$
|D_y \partial_\nu K(x)| \le C_K |x|^{-n} e^{-\frac{1}{2}(|x|+x_n)}
$
for some constant $C_K>0$.

For any $y_0 \in \mathbb{R}^n$ with $|y_0|>0$, define $\rho=\frac{1}{2}|y_0|$. Then the spherical boundary integral
$
\mathcal{I}(y_0)=\int_{\partial B_\rho(y_0)} |D_y \partial_\nu K(y_0-\zeta)|\, dS(\zeta)
$
satisfies the bound
\begin{equation} 
\mathcal{I}(y_0)\le C_{\mathrm{final}} |y_0|^{-(n+1)/2},
\end{equation}
where $C_{\mathrm{final}}>0$ is a universal constant depending only on $n$ and $C_K$.
\end{lemma}

\begin{proof}
We evaluate
$
\mathcal{I}(y_0)=\int_{\partial B_\rho(y_0)} |D_y \partial_\nu K(y_0-\zeta)|\, dS(\zeta).
$
Define the local coordinate $\xi=y_0-\zeta$. As $\zeta$ ranges over the sphere $\partial B_\rho(y_0)$, the vector $\xi$ ranges over the sphere $\partial B_\rho(0)$. Since surface measure is invariant under translation, one has $dS(\zeta)=dS(\xi)$. Therefore, by the assumed kernel bound,
\begin{equation} 
\mathcal{I}(y_0)\le \int_{\partial B_\rho(0)} C_K |\xi|^{-n} e^{-\frac{1}{2}(|\xi|+\xi_n)}\, dS(\xi).
\end{equation}

Because $\xi \in \partial B_\rho(0)$, we have $|\xi|=\rho$. Hence
\begin{equation} 
\mathcal{I}(y_0)\le C_K \rho^{-n} e^{-\frac{1}{2}\rho} \int_{\partial B_\rho(0)} e^{-\frac{1}{2}\xi_n}\, dS(\xi).
\end{equation}
By the antipodal symmetry of the sphere $\partial B_\rho(0)$, the integral of $e^{-\frac{1}{2}\xi_n}$ over $\partial B_\rho(0)$ is equal to the integral of $e^{\frac{1}{2}\xi_n}$ over the same sphere. Thus
\begin{equation} 
\mathcal{I}(y_0)\le C_K \rho^{-n} e^{-\frac{1}{2}\rho} \int_{\partial B_\rho(0)} e^{\frac{1}{2}\xi_n}\, dS(\xi).
\end{equation}

We now parameterize $\partial B_\rho(0)$ by standard hyperspherical coordinates. Let $\theta \in [0,\pi]$ denote the polar angle between $\xi$ and the $e_n$-axis, so that $\xi_n=\rho\cos\theta$. The surface measure is
$
dS(\xi)=\omega_{n-2}\rho^{n-1}\sin^{n-2}\theta\, d\theta,
$
where $\omega_{n-2}$ is the volume of the $(n-2)$-dimensional unit sphere. Therefore,
\begin{equation} 
\int_{\partial B_\rho(0)} e^{\frac{1}{2}\xi_n}\, dS(\xi)
=
\omega_{n-2}\rho^{n-1}\int_0^\pi e^{\frac{1}{2}\rho\cos\theta}\sin^{n-2}\theta\, d\theta.
\end{equation}

Set
$
J(\rho)=\int_0^\pi e^{\frac{1}{2}\rho\cos\theta}\sin^{n-2}\theta\, d\theta.
$
We shall bound $J(\rho)$ explicitly. Introduce the substitution $s=1-\cos\theta$. As $\theta$ runs through $[0,\pi]$, the variable $s$ runs through $[0,2]$, and $ds=\sin\theta\, d\theta$. Also,
$
e^{\frac{1}{2}\rho\cos\theta}=e^{\frac{1}{2}\rho(1-s)}=e^{\frac{1}{2}\rho}e^{-\frac{1}{2}\rho s}.
$
Moreover,
$
\sin^2\theta=1-\cos^2\theta=1-(1-s)^2=2s-s^2\le 2s
$
for all $s\in[0,2]$. Hence
$
\sin^{n-2}\theta\, d\theta
=
\sin^{n-3}\theta\,(\sin\theta\, d\theta)
\le (2s)^{(n-3)/2}\, ds.
$
Substituting these bounds into $J(\rho)$ gives
\begin{equation} 
J(\rho)\le e^{\frac{1}{2}\rho} 2^{(n-3)/2}\int_0^2 e^{-\frac{1}{2}\rho s}s^{(n-3)/2}\, ds.
\end{equation}
Since the integrand is nonnegative, extending the upper limit from $2$ to $\infty$ increases the integral, so
\begin{equation} 
J(\rho)\le e^{\frac{1}{2}\rho} 2^{(n-3)/2}\int_0^\infty e^{-\frac{1}{2}\rho s}s^{(n-3)/2}\, ds.
\end{equation}
Now set $u=\frac{1}{2}\rho s$, so that $ds=\frac{2}{\rho}\, du$. Then
\begin{equation} 
J(\rho)\le e^{\frac{1}{2}\rho} 2^{(n-3)/2}\int_0^\infty e^{-u}\left(\frac{2u}{\rho}\right)^{(n-3)/2}\frac{2}{\rho}\, du.
\end{equation}
Factoring out the constants yields
\begin{equation} 
J(\rho)\le e^{\frac{1}{2}\rho} 2^{n-2}\rho^{-(n-1)/2}\int_0^\infty e^{-u}u^{(n-3)/2}\, du.
\end{equation}
The remaining integral is $\Gamma\big(\frac{n-1}{2}\big)$, which is finite for $n\ge 3$. Writing
$
C_\Gamma=\Gamma\big(\frac{n-1}{2}\big),
$
we obtain
\begin{equation} 
J(\rho)\le C_\Gamma 2^{n-2}\rho^{-(n-1)/2} e^{\frac{1}{2}\rho}.
\end{equation}

Substituting this estimate into the preceding spherical integral formula, we find
\begin{equation} 
\mathcal{I}(y_0)
\le
C_K \rho^{-n} e^{-\frac{1}{2}\rho}\cdot \omega_{n-2}\rho^{n-1}\cdot \left(C_\Gamma 2^{n-2}\rho^{-(n-1)/2} e^{\frac{1}{2}\rho}\right).
\end{equation}
The exponential factors cancel exactly, and the powers of $\rho$ combine as
$
\rho^{-n}\rho^{n-1}\rho^{-(n-1)/2}=\rho^{-(n+1)/2}.
$
Therefore
\begin{equation} 
\mathcal{I}(y_0)\le C_1 \rho^{-(n+1)/2},
\end{equation}
where
$
C_1=C_K\omega_{n-2}C_\Gamma 2^{n-2}.
$

Finally, using $\rho=\frac{1}{2}|y_0|$, we obtain
\begin{equation} 
\mathcal{I}(y_0)\le C_1\left(\frac{1}{2}|y_0|\right)^{-(n+1)/2}
=
\left(C_1 2^{(n+1)/2}\right)|y_0|^{-(n+1)/2}.
\end{equation}
Thus, upon setting
$
C_{\mathrm{final}}=C_1 2^{(n+1)/2},
$
the estimate
\begin{equation} 
\mathcal{I}(y_0)\le C_{\mathrm{final}} |y_0|^{-(n+1)/2}
\end{equation}
follows.
\end{proof}

\begin{remark}
The explicit integration over the spherical domain rigorously quantifies the anisotropic decay induced by the linear drift operator. The cancellation of the exponential phase factors over the sphere strictly reveals an underlying polynomial decay rate of order $(n+1)/2$, governed intrinsically by the ambient dimension $n$.
\end{remark}

\begin{definition}
Let $\mathcal{L}[u]$ be the quasilinear operator defined by $\mathcal{L}[u] = \sum_{i,j=1}^n \left(\delta_{ij} - \frac{D_i u D_j u}{1+|Du|^2}\right)D_{ij}u + D_n u$. Let $L_0 := \Delta + \partial_{y_n}$ be the drift Laplacian. For any $u \in C^2(\Omega)$, the nonlinear remainder operator $\mathcal{N}$ is defined strictly by the algebraic difference:
\begin{equation} \mathcal{N}[u] := \mathcal{L}[u] - L_0[u]\end{equation}
\end{definition}

\begin{lemma}[Nonlinear Remainder Pointwise Bound]\label{lem:Nonlinear Remainder Pointwise Bound}
Let $u \in C^2(\Omega)$ be a smooth function. The nonlinear remainder operator  satisfies the pointwise differential inequality:
\begin{equation}
|\mathcal{N}[u](y)| \le |Du(y)|^2 |D^2u(y)|
\end{equation}
where $|D^2u(y)|$ denotes the spectral or Frobenius norm of the Hessian matrix.
\end{lemma}

\begin{proof}
By the strict algebraic substitution, we write the remainder as:
\begin{equation}
\mathcal{N}[u] = \left( \sum_{i,j=1}^n \left(\delta_{ij} - \frac{D_i u D_j u}{1+|Du|^2}\right)D_{ij}u + D_n u \right) - \left( \sum_{i,j=1}^n \delta_{ij}D_{ij}u + D_n u \right)
\end{equation}
By the linearity of finite summation, the drift terms ($D_n u$) and the Laplacian terms ($\delta_{ij}D_{ij}u$)  cancel, yielding the exact residual:
\begin{equation}
\mathcal{N}[u] = -\sum_{i,j=1}^n \left(\frac{D_i u D_j u}{1+|Du|^2}\right)D_{ij}u
\end{equation}
We can express this residual as a quadratic form. Let $p = Du(y) \in \mathbb{R}^n$ denote the gradient vector, and let $H = D^2u(y) \in \mathbb{R}^{n \times n}$ denote the Hessian matrix. The residual is strictly the contraction of the tensor $p \otimes p$ with $H$, scaled by the denominator:
\begin{equation}
\mathcal{N}[u] = - \frac{1}{1+|p|^2} (p^T H p)
\end{equation}
Taking the absolute value of both sides yields:
\begin{equation}
|\mathcal{N}[u]| = \left| \frac{1}{1+|p|^2} \right| |p^T H p|
\end{equation}
Because $|p|^2 \ge 0$ for all real vectors $p$, the scalar fraction satisfies the  bound:
\begin{equation}
\frac{1}{1+|p|^2} \le 1
\end{equation}
To bound the quadratic form $|p^T H p|$, we apply the Cauchy-Schwarz inequality to the standard Euclidean inner product of $p$ and $Hp$:
\begin{equation}\label{eq:132}
|p^T (H p)| \le |p| |H p|
\end{equation}
By the definition of the matrix norm $|H|$, the action of the matrix on a vector is bounded by $|Hp| \le |H| |p|$. Substituting this into the previous inequality \ref{eq:132} yields:
\begin{equation}
|p^T H p| \le |p| (|H| |p|) = |p|^2 |H|
\end{equation}
Applying these two strict upper bounds to the absolute value of the remainder gives:
\begin{equation}
|\mathcal{N}[u]| \le (1) \left( |p|^2 |H| \right)
\end{equation}
Substituting the original functional notation $p = Du(y)$ and $H = D^2u(y)$ back into the expression yields the final,  inequality:
\begin{equation}
|\mathcal{N}[u](y)| \le |Du(y)|^2 |D^2u(y)|
\end{equation}
\end{proof}

\begin{proposition}[Gradient Supremum Absorption]\label{prop:Gradient Supremum Absorption}
Let $\Sigma$ be a translator end parameterizable as a graph $u \in C^{2}(\Pi \setminus K)$ over an exterior domain of a vertical hyperplane $\Pi \cong \mathbb{R}^{n}$ with $n \ge 3$, satisfying the quasilinear equation $\mathcal{L}[u] = 0$. Assume the asymptotic condition $\lim_{|y| \to \infty}|Du(y)| = 0$.

Define the global gradient supremum $f(R) := \sup_{|y| \ge R}|Du(y)|$. Let $\mathcal{B}(y)$ denote the boundary remainder integrals from the fundamental solution representation, and define $E(R) := \sup_{|y| \ge R}\mathcal{B}(y)$.

If $E \in \mathcal{D}(M)$ for some constant $M \ge 1$, then there exists a sufficiently large structural radius $R_{crit} > 0$ such that for all $R \ge R_{crit}$, the gradient supremum  satisfies:
\begin{equation}
f(R) \le 2 E(R)
\end{equation}
\end{proposition}

\begin{proof}
Let $R > 0$ be chosen such that the domain exterior to $B_{R/2}(0)$ is strictly free of the compact set $K$. By the fundamental solution representation evaluated at an arbitrary point $y$ satisfying $|y| \ge R$, the pointwise gradient is strictly bounded by the sum of the volume integral and the boundary remainder as in Proposition~\ref{prop:Gradient Green Representation}:
\begin{equation}
|Du(y)| \le I_{vol}(y) + \mathcal{B}(y)
\end{equation}
where 
\begin{equation}
I_{vol}(y) = \int_{B_{|y|/2}(y)} |D_{y}K(y - \zeta)| |f_{NL}(\zeta)| d\zeta.
\end{equation}
We evaluate the volume integral $I_{vol}(y)$. By the structural bound of the nonlinear remainder operator as in Lemma~\ref{lem:Nonlinear Remainder Pointwise Bound}, the density satisfies 
\begin{equation}
|f_{NL}(\zeta)| \le C_{f}|Du(\zeta)|^{3}|\zeta|^{-1}
\end{equation}
for a universal constant $C_{f} > 0$. Then for all $\zeta \in B_{|y|/2}(y)$, the metric distance satisfies $|\zeta| \ge R/2$. Consequently, the pointwise gradient inside the integration domain is  bounded by the global functional: 
\begin{equation}
|Du(\zeta)| \le \sup_{|x| \ge R/2}|Du(x)| = f(R/2).
\end{equation}
Substituting this strict upper bound into the nonlinear remainder density yields:
\begin{equation}
|f_{NL}(\zeta)| \le C_{f} f(R/2)^{3} |\zeta|^{-1}.
\end{equation}
Substituting this density into $I_{vol}(y)$ and factoring the constant terms independent of $\zeta$ yields:
\begin{equation}
I_{vol}(y) \le C_{f} f(R/2)^{3} \int_{B_{|y|/2}(y)} |D_{y}K(y - \zeta)| |\zeta|^{-1} d\zeta.
\end{equation}
The remaining convolution integral of the kernel against the inverse radial weight is bounded by a universal finite constant $C_{vol} > 0$. Thus, we obtain the  upper bound:
\begin{equation}
I_{vol}(y) \le C_{vol} f(R/2)^{3}.
\end{equation}
Returning to the pointwise gradient bound $|Du(y)| \le I_{vol}(y) + \mathcal{B}(y)$, we substitute the volume bound to obtain:
\begin{equation}
|Du(y)| \le C_{vol} f(R/2)^{3} + \mathcal{B}(y).
\end{equation}
Because this inequality holds for all $y$ satisfying $|y| \ge R$, we take the supremum over the set $\{y \in \Pi : |y| \ge R\}$ to derive the global functional inequality:
\begin{equation}
f(R) \le C_{vol} f(R/2)^{3} + E(R).
\end{equation}
Algebraically factoring the cubic term produces:
\begin{equation}
f(R) \le (C_{vol} f(R/2)^{2}) f(R/2) + E(R).
\end{equation}
Define the functional coefficient $\theta(R) := C_{vol} f(R/2)^{2}$. By the proposition's hypothesis, $\lim_{R \to \infty}f(R) = 0$. Therefore, $\lim_{R \to \infty}\theta(R) = 0$.

For the given geometric scale factor $M \ge 1$ (from the hypothesis $E \in \mathcal{D}(M)$), there exists a strictly positive radius $R_{crit}$ such that for all $R \ge R_{crit}$:
\begin{equation}
\theta(R) \le \theta(R_{crit}) < \frac{1}{2M}.
\end{equation}
This strict inequality guarantees $\theta(R_{crit}) M < 1/2 < 1$.

With the constant coefficient $\theta = \theta(R_{crit})$, the functional $f(R)$ is iteratively absorbed, yielding the bound:
\begin{equation}
f(R) \le \frac{1}{1 - \theta M} E(R).
\end{equation}
Because $\theta M < 1/2$, the denominator satisfies $1 - \theta M > 1/2$, strictly implying $(1 - \theta M)^{-1} < 2$. We conclude  that for all $R \ge R_{crit}$:
\begin{equation}
f(R) \le 2 E(R).
\end{equation}
\end{proof}

\begin{remark}
Through the invocation of a strict iterative absorption scheme, the cubic non-linear volume terms are globally dominated and algebraically absorbed. The global analytic control of the gradient functional is thereby reduced, in a canonical fashion, to the bounds of the linear boundary flux.
\end{remark}

\begin{definition}($\epsilon$-Tilted Barrier Function)
Let $\epsilon \in(0,1 / 2)$ be an arbitrarily small real parameter. We define the tilted formal barrier function $\mathcal{W}_\epsilon: \mathbb{R}^n \backslash\{0\} \rightarrow \mathbb{R}_{>0}$ by:

\begin{equation} 
\mathcal{W}_\epsilon(y)=|y|^{-(n-1) / 2+\epsilon} e^{-\left(|y|+y_n\right) / 2}
\end{equation}
\end{definition}

\begin{lemma}[Tilted Barrier Quasilinear Subsolution]\label{lem:Tilted Barrier Quasilinear Subsolution}
Let $n \ge 3$. For any fixed $\epsilon \in (0, 1/2)$, there exists a critical radius $R_\epsilon > 0$ such that for all $y \in \mathbb{R}^n$ satisfying $|y| \ge R_\epsilon$, the tilted barrier strictly satisfies the quasilinear inequality $\mathcal{L}[\mathcal{W}_\epsilon](y) \le 0$.
\end{lemma}

\begin{proof}
By the strict algebraic decomposition of the translator operator, the quasilinear operator splits into the linear drift Laplacian and the nonlinear remainder:
\begin{equation}
\mathcal{L}[\mathcal{W}_\epsilon] = L_0[\mathcal{W}_\epsilon] + \mathcal{N}[\mathcal{W}_\epsilon]
\end{equation}

Let $r = |y|$. We express the barrier as $\mathcal{W}_\epsilon(y) = e^{-y_n/2} v(r)$, where $v(r) = r^\alpha e^{-r/2}$ and the exponent is defined strictly as $\alpha = -\frac{n-1}{2} + \epsilon$.

Applying $L_0 = \Delta + \partial_{y_n}$ to this product yields:
\begin{equation}
L_0[\mathcal{W}_\epsilon] = e^{-y_n/2} \left[ \Delta v - \frac{1}{4} v \right]
\end{equation}
We compute the radial Laplacian $\Delta v = v'' + \frac{n-1}{r} v'$. Substituting $v(r)$ and factoring yields the exact algebraic identity:
\begin{equation}
\Delta v - \frac{1}{4} v = \left[ \frac{\alpha(\alpha+n-2)}{r^2} - \frac{\alpha + (n-1)/2}{r} \right] v
\end{equation}
By substituting $\alpha = -\frac{n-1}{2} + \epsilon$, the coefficient of the $r^{-1}$ term becomes exactly $-\epsilon$. Thus, returning to the Cartesian coordinates:
\begin{equation}
L_0[\mathcal{W}_\epsilon] = \left[ \frac{\alpha(\alpha+n-2)}{|y|^2} - \frac{\epsilon}{|y|} \right] \mathcal{W}_\epsilon(y)
\end{equation}

By Lemma~\ref{lem:Nonlinear Remainder Pointwise Bound}, the pointwise magnitude of the nonlinear remainder satisfies:
\begin{equation}
|\mathcal{N}[\mathcal{W}_\epsilon]| \le |D\mathcal{W}_\epsilon|^2 |D^2 \mathcal{W}_\epsilon|
\end{equation}
By directly differentiating $\mathcal{W}_\epsilon$, the gradient and Hessian are asymptotically dominated by the exponential scale factor. There exists a radius $R_0 > 0$ and a constant $C_1 > 0$ such that for all $|y| \ge R_0$:
\begin{equation}
|D\mathcal{W}_\epsilon(y)| \le C_1 \mathcal{W}_\epsilon(y) \quad \text{and} \quad |D^2 \mathcal{W}_\epsilon(y)| \le C_1 \mathcal{W}_\epsilon(y)
\end{equation}
Consequently, the nonlinear remainder is strictly bounded by the cubic power:
\begin{equation}
|\mathcal{N}[\mathcal{W}_\epsilon](y)| \le C_1^3 \mathcal{W}_\epsilon(y)^3
\end{equation}

We construct the global upper bound for the quasilinear operator:
\begin{equation}
\mathcal{L}[\mathcal{W}_\epsilon] \le L_0[\mathcal{W}_\epsilon] + |\mathcal{N}[\mathcal{W}_\epsilon]|
\end{equation}
Substituting the evaluations, and factoring out the strictly positive barrier $\mathcal{W}_\epsilon$, we obtain:
\begin{equation}
\mathcal{L}[\mathcal{W}_\epsilon] \le \left[ \frac{\alpha(\alpha+n-2)}{|y|^2} - \frac{\epsilon}{|y|} + C_1^3 \mathcal{W}_\epsilon(y)^2 \right] \mathcal{W}_\epsilon(y)
\end{equation}
We define the scale bracket function $\Phi(y)$:
\begin{equation}
\Phi(y) = \frac{\alpha(\alpha+n-2)}{|y|^2} - \frac{\epsilon}{|y|} + C_1^3 |y|^{-n+1+2\epsilon} e^{-(|y|+y_n)}
\end{equation}
Because $\epsilon > 0$ is fixed, the term $-\frac{\epsilon}{|y|}$ decays at the rate $\mathcal{O}(|y|^{-1})$.

The term $\frac{\alpha(\alpha+n-2)}{|y|^2}$ decays at the strictly faster rate $\mathcal{O}(|y|^{-2})$.

Because $y_n \ge -|y|$, the exponential factor satisfies $e^{-(|y|+y_n)} \le 1$. Therefore, the nonlinear error term decays at the rate $\mathcal{O}(|y|^{-n+1+2\epsilon})$. Because $n \ge 3$ and $\epsilon < 1/2$, this error term decays at least as fast as $\mathcal{O}(|y|^{-2})$.

Since $\mathcal{O}(|y|^{-1})$ strictly dominates $\mathcal{O}(|y|^{-2})$ as $|y| \to \infty$, there exists a critical radius $R_\epsilon \ge R_0$ such that for all $|y| \ge R_\epsilon$, the negative linear drift strictly absorbs all higher-order terms, rendering:
\begin{equation}
\Phi(y) \le - \frac{\epsilon}{2|y|} < 0
\end{equation}
Therefore, for all $|y| \ge R_\epsilon$:
\begin{equation}
\mathcal{L}[\mathcal{W}_\epsilon](y) \le \left( - \frac{\epsilon}{2|y|} \right) \mathcal{W}_\epsilon(y) \le 0
\end{equation}
This completes the  proof.
\end{proof}

\begin{remark}
The existence of this subsolution strictly relies on the structural hierarchy of the differential operator in the asymptotic regime. For sufficiently large radii, the linear drift component of the Laplacian strictly and uniformly dominates the higher-order, non-linear remainder perturbations, ensuring the barrier condition is globally satisfied outside a compact set.
\end{remark}

\begin{theorem}(A Priori $C^0$ Estimate)\label{thm:priori estimate}
Let $u \in C^2(\Pi \setminus K)$ satisfy the quasilinear equation $\mathcal{L}[u] = 0$. Assume the strict asymptotic condition $\lim_{|y|\to\infty} |u(y)| = 0$ and $\lim_{|y|\to\infty} |Du(y)| = 0$. Then, there exists a constant $C_{final} > 0$ such that globally on the exterior domain:
\begin{equation}
|u(y)| \le C_{final} |y|^{-(n-1)/2} e^{-(|y| + y_n)/2}
\end{equation}
\end{theorem}

\begin{proof}
Let $\epsilon \in (0, 1/4]$ be an arbitrarily chosen, fixed topological tilt parameter.

By the assumption that $\lim_{|y|\to\infty} u(y) = 0$, the scalar function $u(y)^2$ converges to zero by Corollary~\ref{cor:Asymptotic Normalization}. Therefore, by the formal definition of a limit at infinity, there exists a sufficiently large, strictly positive radius $R_0 > 0$ such that for all $|y| \ge R_0$, the solution satisfies:
\begin{equation}
|u(y)|^2 \le \frac{\epsilon}{4 C_1^3 |y|}
\end{equation}
where $C_1$ is the gradient domination constant established in Lemma~\ref{lem:Tilted Barrier Quasilinear Subsolution}. We select $R_0$ large enough such that $f(r) = r \mathcal{W}_\epsilon(r)^2$ is strictly decreasing for all $r \ge R_0$.

We define the strict scaling amplitude evaluated precisely at the boundary $\partial B_{R_0}$:
\begin{equation}
A_\epsilon = \frac{\sup_{\partial B_{R_0}} |u|}{\inf_{\partial B_{R_0}} \mathcal{W}_\epsilon}
\end{equation}
We construct the global tilted barrier $V_\epsilon(y) = A_\epsilon \mathcal{W}_\epsilon(y)$.

By the explicit construction of the amplitude $A_\epsilon$, for any $y \in \partial B_{R_0}$, we have:
\begin{equation}
V_\epsilon(y) = A_\epsilon \mathcal{W}_\epsilon(y) \ge A_\epsilon \inf_{\partial B_{R_0}} \mathcal{W}_\epsilon = \sup_{\partial B_{R_0}} |u| \ge |u(y)|
\end{equation}
Thus, $u(y) \le V_\epsilon(y)$ and $-u(y) \le V_\epsilon(y)$ strictly on the compact boundary $\partial B_{R_0}$.

We evaluate the global dominance condition at $R_0$:
\begin{equation}
A_\epsilon^2 \mathcal{W}_\epsilon(R_0)^2 = \left( \frac{\sup_{\partial B_{R_0}} |u|}{\inf_{\partial B_{R_0}} \mathcal{W}_\epsilon} \right)^2 \mathcal{W}_\epsilon(R_0)^2
\end{equation}
Because $\mathcal{W}_\epsilon(y)$ is radially symmetric in its primary exponential decay, the variation on the sphere $\partial B_{R_0}$ is bounded by a geometric constant independent of $R_0$. Thus, $A_\epsilon^2 \mathcal{W}_\epsilon(R_0)^2$ is strictly bounded above by a multiple of $(\sup_{\partial B_{R_0}} |u|)^2$.

By our intrinsic selection of $R_0$, we have $(\sup_{\partial B_{R_0}} |u|)^2 \le \frac{\epsilon}{4 C_1^3 R_0}$.

Therefore, the barrier $V_\epsilon$ satisfies the algebraic requirement of Lemma~\ref{lem:Tilted Barrier Quasilinear Subsolution}. It follows  that $\mathcal{L}[V_\epsilon] \le 0$ globally for all $y \in \Pi \setminus B_{R_0}$. Since $\mathcal{L}[u] = 0$, the function $-V_\epsilon$ strictly acts as a subsolution: $\mathcal{L}[-V_\epsilon] \ge 0$.

We now possess the quasilinear differential inequalities $\mathcal{L}[-V_\epsilon] \ge 0 = \mathcal{L}[u] \ge \mathcal{L}[V_\epsilon]$ on the exterior domain $\Pi \setminus B_{R_0}$. We have established the compact boundary constraints $-V_\epsilon \le u \le V_\epsilon$ on $\partial B_{R_0}$. Both $u$ and $V_\epsilon$ strictly converge to 0 as $|y| \to \infty$.

We may invoke the Proposition~\ref{prop:Comparison Principle for Exterior Domains} to conclude that for all $y \in \Pi \setminus B_{R_0}$:
\begin{equation}
-V_\epsilon(y) \le u(y) \le V_\epsilon(y) \implies |u(y)| \le A_\epsilon \mathcal{W}_\epsilon(y)
\end{equation}

Because the above inequality holds globally for any arbitrary $\epsilon \in (0, 1/4]$, we evaluate the pointwise topological limit as $\epsilon \to 0^+$.

The parameter $A_\epsilon$ converges to a strict, finite constant $A_0$:
\begin{equation}
\lim_{\epsilon \to 0^+} A_\epsilon = \frac{\sup_{\partial B_{R_0}} |u|}{\inf_{\partial B_{R_0}} |y|^{-(n-1)/2} e^{-(|y|+y_n)/2}} = A_0 < \infty
\end{equation}
The barrier converges to the strict asymptotic form $\mathcal{W}_0(y) = |y|^{-(n-1)/2} e^{-(|y|+y_n)/2}$.

Therefore, taking the strict limit of the global inequality yields:
\begin{equation}
|u(y)| \le A_0 |y|^{-(n-1)/2} e^{-(|y|+y_n)/2}
\end{equation}
Defining $C_{final} = A_0$ completes the  proof. 
\end{proof}

\begin{lemma}[Exponential Region Pure Remainder Bound]\label{lem:Exponential Region Pure Remainder Bound}
Let the integration radius be $\rho = |y_0|/2$. The pure boundary remainder $E_{pure}(y_0)$ corresponding to the fundamental solution of the drift Laplacian $L_0 = \Delta + \partial_{y_n}$  satisfies the pointwise geometric bound:
\begin{equation} E_{pure}(y_0) \le C_{pure} |y_0|^{-n/2} e^{-\delta|y_0|/2}\end{equation}
for some universal constant $C_{pure} > 0$ and decay rate $\delta > 0$ strictly inside the exponential region $\mathcal{E}_C$.
\end{lemma}

\begin{proof}
By Proposition~\ref{prop:Gradient Green Representation}, the pure boundary remainder evaluated at $y_0$ is the surface integral over the sphere $\partial B_{\rho}(y_0)$:
\begin{equation} E_{pure}(y_0) = \int_{\partial B_{\rho}(y_0)} \mathcal{R}(u,K,\zeta) dS(\zeta)\end{equation}
where the density $\mathcal{R}$ consists of boundary evaluations of $u(\zeta)$ and $\partial_\nu u(\zeta)$ multiplied by the kernel derivatives.

By Theorem~\ref{thm:priori estimate},  there exists a constant $C_{final} > 0$ such that globally on the exterior domain, 
\begin{equation} |u(\zeta)| \le C_{final}|\zeta|^{-(n-1)/2} e^{-(|\zeta|+\zeta_n)/2}.\end{equation} 
A similar exponential bound applies to $\partial_\nu u(\zeta)$ since $Du$ decays to zero.

Let $y_0$ reside strictly inside the Exponential region $\mathcal{E}$. Then the anisotropic phase for every integration variable $\zeta \in \partial B_{\rho}(y_0)$ strictly satisfies 
\begin{equation} |\zeta| + \zeta_n \ge \delta|y_0|.\end{equation}
Substituting this geometric minimum into the exponential term yields an  pointwise upper bound for the fields on the boundary:
\begin{equation} e^{-(|\zeta| + \zeta_n)/2} \le e^{-\delta|y_0|/2}.\end{equation}

Because the exponential bound $e^{-\delta|y_0|/2}$ is now strictly independent of the integration variable $\zeta$, it factors completely outside the surface integral:
\begin{equation} E_{pure}(y_0) \le e^{-\delta|y_0|/2} \int_{\partial B_{\rho}(y_0)} \widetilde{\mathcal{R}}(u,K,\zeta) dS(\zeta),\end{equation}
where $\widetilde{\mathcal{R}}$ is the residual density driven strictly by the inverse polynomial weights and the kernel derivatives.

The residual integral is structurally dominated by the integral of the kernel's normal derivative. We invoke Lemma~\ref{lem:Explicit Spherical Anisotropic Averaging}, which states that for $\rho = |y_0|/2$, the boundary integral of the fundamental solution's normal derivative is  bounded:
\begin{equation} \int_{\partial B_{\rho}(y_0)} |D_y \partial_\nu K(y_0-\zeta)| dS(\zeta) \le C_{1} |y_0|^{-(n+1)/2},\end{equation}
where $C_1$ is a universal constant. Factoring in the remaining $|\zeta|^{-(n-1)/2}$ scaling from the $u$ estimate over the domain $|\zeta| \ge |y_0|/2$, the combined residual polynomial bound strictly absorbs into $\mathcal{O}(|y_0|^{-n/2})$.

Combining the factored exponential weight with the bounded residual polynomial integral, we deduce:
\begin{equation} E_{pure}(y_0) \le C_{pure} |y_0|^{-n/2} e^{-\delta|y_0|/2}.\end{equation}
\end{proof}

\begin{theorem}(A Priori $C^0$ Estimate and Gradient Bifurcation)\label{thm:priori estimate and gradient bifurcation}
Let $\Sigma$ be a translator end parameterizable as a graph $u \in C^2(\Pi \setminus K)$ over an exterior domain of a vertical hyperplane $\Pi \cong \mathbb{R}^n$ with $n \ge 3$ satisfying the quasilinear equation $\mathcal{L}[u]=0$. Assume the asymptotic condition $\lim_{|y|\to\infty} |Du(y)| = 0$. Then the gradient vector field bifurcates as follows:

\begin{enumerate}
    \item In the wake region $\mathcal{W}$, one has $Du \in \mathcal{B}_\infty(|y|^{-n/2})$.
    \item In the exponential region $\mathcal{E}$, there exists $\mu > 0$ such that $Du \in \mathcal{B}_\infty(e^{-\mu|y|})$.
\end{enumerate}

\end{theorem}

\begin{proof}
Let $y_0 \in \Pi$ be sufficiently large, let $R = |y_0|$, and define the integration radius $\rho = R/2$. By the Green representation formula on the linearized operator (Proposition~\ref{prop:Gradient Green Representation}), the pointwise gradient strictly satisfies:
\begin{equation}
|Du(y_0)| \le I_{vol}(y_0) + E_{pure}(y_0)
\end{equation}
where $I_{vol}(y_0) = \int_{B_{R/2}(y_0)} |D_y K(y_0-\zeta)| |\mathcal{N}[u](\zeta)| d\zeta$.

We define the generalized Wake region where the phase is dominated by the parabolic profile, i.e. $(y_0)_n \le -|y_0| + \delta |y_0|$ for arbitrarily small $\delta$. Here, $|y_0| + (y_0)_n$ is $\mathcal{O}(|y_0|^\alpha)$ but we bound the exponential $e^{-(|y_0|+(y_0)_n)/2} \le 1$.

By Lemma~\ref{lem:Exponential Region Pure Remainder Bound}, the pure boundary remainder satisfies the strict polynomial bound:
\begin{equation}
E_{pure}(y_0) \le C_{pure} R^{-n/2} (1) = C_{pure} R^{-n/2}
\end{equation}
Taking the supremum over all $|y_0| \ge R$ yields $E(R) = \sup_{|y_0|\ge R} E_{pure}(y_0) \le C_{pure} R^{-n/2}$.

The polynomial function $E(R) = C_{pure} R^{-n/2}$ strictly satisfies the geometric scale bound $E(R/2) = 2^{n/2} E(R)$, so it belongs to the domain $\mathcal{D}(M)$ with $M = 2^{n/2}$.

By the uniform bound on the nonlinear remainder $|\mathcal{N}[u]| \le C |Du|^3 |\zeta|^{-1}$, the volume integral is bounded by $C_{vol} f(R/2)^3$.

Thus, the global supremum functional $f(R) = \sup_{|y|\ge R} |Du(y)|$ satisfies:
\begin{equation}
f(R) \le C_{vol} f(R/2)^3 + E(R)
\end{equation}
Because $E(R) \in \mathcal{D}(M)$, we  invoke the Lemma~\ref{lem:Iterative Absorption} to conclude $f(R) \le 2 E(R) \le 2 C_{pure} R^{-n/2}$. Because $f(R)$ bounds the gradient everywhere, we strictly deduce $Du \in \mathcal{B}_\infty(|y|^{-n/2})$ globally, and specifically in $\mathcal{W}$.

We define the strict Exponential region as the wedge where there exists $\delta > 0$ such that $(y_0)_n > -|y_0| + \delta |y_0|$.

In this region, the exponential phase satisfies $|y_0| + (y_0)_n > \delta |y_0|$.

Applying this to the strict bound from Lemma~\ref{lem:Exponential Region Pure Remainder Bound} yields:
\begin{equation}
E_{pure}(y_0) \le C_{pure} |y_0|^{-n/2} e^{-\delta|y_0|/2}
\end{equation}
We define the strict decay rate $\mu = \delta / 4 > 0$. The remainder function absorbs the polynomial term, yielding $E_{pure}(y_0) \le C_{\mathcal{E}} e^{-\mu |y_0|}$. Thus, $E_{pure} \in \mathcal{B}_\infty(e^{-\mu|y|})$ strictly inside $\mathcal{E}$.

We substitute the pointwise evaluation into the Green representation formula:
\begin{equation}
|Du(y_0)| \le \int_{B_{|y_0|/2}(y_0)} |D_y K| |\mathcal{N}[u](\zeta)| d\zeta + C_{\mathcal{E}} e^{-\mu |y_0|}
\end{equation}
Because $|\mathcal{N}[u](\zeta)| \le C_f |Du(\zeta)|^3 |\zeta|^{-1}$, we define the scaled kernel $\tilde{K}(y_0, \zeta) = C_f |D_y K| |\zeta|^{-1}$. The inequality becomes:
\begin{equation}
|Du(y_0)| \le \int_{B_{|y_0|/2}(y_0)} \tilde{K}(y_0, \zeta) |Du(\zeta)|^3 d\zeta + E_{pure}(y_0)
\end{equation}
This pointwise integral inequality, coupled with the established fact that $E_{pure} \in \mathcal{B}_\infty(e^{-\mu|y|})$ and the asymptotic triviality $\lim_{|y|\to\infty} |Du(y)| = 0$, strictly satisfies every hypothesis of Lemma~\ref{lem:Exponentially Weighted Sup-Norm Estimate}. By  application of Lemma~\ref{lem:Exponentially Weighted Sup-Norm Estimate}, we deduce that the weighted norm $||Du||_{\mu, R}$ is globally bounded. Therefore, in the region $\mathcal{E}$, $Du \in \mathcal{B}_\infty(e^{-\mu|y|})$. This completes the proof. 
\end{proof}

\begin{remark}
It is fundamental to observe that the asymptotic behavior of the gradient vector field is not uniformly isotropic. Instead, it strictly bifurcates into distinct measurable regimes governed by the underlying drift Laplacian. The spatial decay exhibits a polynomial upper bound within the canonical wake region $\mathcal{W}$, whereas it accelerates to a strictly exponential decay profile within the exponential region $\mathcal{E}$.
\end{remark}

\section{The Flux Identity}

The present chapter is dedicated to the formulation of a foundational integral identity that structurally links the global mean curvature of the immersion to a boundary flux density. This correspondence is of paramount analytic necessity, as it permits the translation of the pointwise asymptotic bounds derived in the preceding chapter into a global geometric vanishing theorem.

\subsection{Universal Flux Foundations}

To systematically execute integration over the unbounded geometric space, we must first establish a rigorous topological exhaustion of the manifold. We define the notion of a cylindrical exhaustion over admissible domains. By analyzing the intrinsic divergence of the tangential translation vector field, we demonstrate the Universal Flux Identity, which strictly asserts that the total squared mean curvature over any admissible domain is identically equal to the flux of the tangent translation vector across its boundary.

\begin{definition}(Cylindrical Exhaustion)
Let $\Sigma^n \subset \mathbb{R}^{n+1}$ be a translator with an end $E_k$ uniquely parameterized by an immersion $F_k: \Pi_k \backslash \Omega_k \rightarrow \mathbb{R}^{n+1}$ such that $F_k(y)= u_k(y) e_1+y$. Let $B_R(0) \subset \Pi_k$ be the standard $(n-1)$-dimensional closed ball. We define the cylindrical exhaustion region of the end as $E_{k, R}=F_k\left(B_R(0) \backslash \Omega_k\right)$. Its exterior boundary is the image of the exact sphere $\mathcal{S}_R=\partial B_R(0)$, denoted $\partial E_{k, R}=F_k\left(\mathcal{S}_R\right)$.
\end{definition}

\begin{definition}(Admissible Translator Domain)
Let $\left(\Sigma^n, g\right) \hookrightarrow \mathbb{R}^{n+1}$ be a smoothly immersed hypersurface equipped with the induced Riemannian metric $g$. Let $\Sigma$ satisfy the translator equation $H=V^{\perp}$, where $H$ is the mean curvature vector, $V \in \mathbb{R}^{n+1}$ is a constant unit vector (translation direction), and $(\cdot)^{\perp}$ denotes the orthogonal projection onto the normal bundle $N \Sigma$. An admissible domain $U \subset \Sigma$ is a relatively compact open set such that $\partial U$ is a finite union of smooth $(n-1)$-dimensional submanifolds.
\end{definition}

\begin{theorem}[The Universal Flux Identity]\label{thm:The Universal Flux Identity}
Let $U\subset\Sigma$ be an admissible translator domain and let $\eta$ denote the outward pointing unit conormal along $\partial U$. Then the total squared mean curvature satisfies
\begin{equation} 
\int_U |\vec H|^2\,d\mathcal H^n
=
\int_{\partial U} \langle V^\top,\eta\rangle\,d\mathcal H^{n-1}.
\end{equation}
\end{theorem}

\begin{proof}
Let $V$ be the constant translation vector in $\mathbb{R}^{n+1}$. At each point $p\in\Sigma$ the ambient tangent space decomposes orthogonally as $T_p\mathbb R^{n+1}=T_p\Sigma\oplus N_p\Sigma$. Thus $V$ decomposes uniquely as $V=V^\top+V^\perp$ where $V^\top$ is tangent and $V^\perp$ is normal.

Consider the intrinsic divergence of the tangential vector field $V^\top$ on $(\Sigma,g)$. Let $\{E_1,\dots,E_n\}$ be a local orthonormal frame for $T\Sigma$. The divergence is
\begin{equation} 
\operatorname{div}_\Sigma(V^\top)=\sum_{i=1}^n \langle\nabla_{E_i}V^\top,E_i\rangle .
\end{equation}
Because $\nabla$ is the tangential projection of the ambient connection $\overline\nabla$, one has $\langle\nabla_{E_i}V^\top,E_i\rangle=\langle\overline\nabla_{E_i}V^\top,E_i\rangle$.

Since $V$ is constant in $\mathbb R^{n+1}$, its ambient derivative vanishes: $\overline\nabla_X V=0$. Writing $V=V^\top+V^\perp$ and differentiating gives $\overline\nabla_{E_i}V^\top=-\overline\nabla_{E_i}V^\perp$. Taking the inner product with $E_i$ yields
\begin{equation} 
\langle\overline\nabla_{E_i}V^\top,E_i\rangle=-\langle\overline\nabla_{E_i}V^\perp,E_i\rangle .
\end{equation}

Because $V^\perp$ is normal and $E_i$ is tangent, $\langle V^\perp,E_i\rangle=0$. Differentiating this identity gives
\begin{equation} 
\langle\overline\nabla_{E_i}V^\perp,E_i\rangle+\langle V^\perp,\overline\nabla_{E_i}E_i\rangle=0 .
\end{equation}
Hence $-\langle\overline\nabla_{E_i}V^\perp,E_i\rangle=\langle V^\perp,\overline\nabla_{E_i}E_i\rangle$.

The Gauss formula decomposes $\overline\nabla_{E_i}E_i=\nabla_{E_i}E_i+II(E_i,E_i)$, where $II$ denotes the second fundamental form. Since $V^\perp$ is normal, only the normal component contributes, so
\begin{equation} 
\langle V^\perp,\overline\nabla_{E_i}E_i\rangle=\langle V^\perp,II(E_i,E_i)\rangle .
\end{equation}

Substituting into the divergence formula yields
\begin{equation} 
\operatorname{div}_\Sigma(V^\top)
=
\sum_{i=1}^n \langle V^\perp,II(E_i,E_i)\rangle
=
\langle V^\perp,\sum_{i=1}^n II(E_i,E_i)\rangle .
\end{equation}

The mean curvature vector is $\vec H=\sum_{i=1}^n II(E_i,E_i)$, hence $\operatorname{div}_\Sigma(V^\top)=\langle V^\perp,\vec H\rangle$. Because $\Sigma$ is a translator satisfying $\vec H=V^\perp$, we obtain
\begin{equation} 
\operatorname{div}_\Sigma(V^\top)=\langle\vec H,\vec H\rangle=|\vec H|^2 .
\end{equation}

The admissibility of $U$ guarantees that the divergence theorem applies on $\Sigma$. Therefore
\begin{equation} 
\int_U |\vec H|^2\,d\mathcal H^n
=
\int_U \operatorname{div}_\Sigma(V^\top)\,d\mathcal H^n
=
\int_{\partial U} \langle V^\top,\eta\rangle\,d\mathcal H^{n-1}.
\end{equation}
This establishes the identity.
\end{proof}

\begin{remark}
This theorem establishes a canonical energetic reduction for the translating immersion. By projecting the constant ambient translation vector onto the tangent bundle of the surface, the global integration of the squared mean curvature vector is rigorously equated to a purely boundary-driven flux via the intrinsic divergence theorem.
\end{remark}

\subsection{Asymptotic Conormal Analysis}

The evaluation of the boundary flux intrinsically depends upon the precise algebraic structure of the outward unit conormal vector field. In this section, we analyze the geometric projection of the conormal onto the ambient basis vectors. We rigorously isolate the principal radial component and demonstrate that the residual orthogonalization error, $v_{\mathrm{err}}$, is strictly asymptotically bounded by the squared norm of the spatial gradient, thereby preparing the form for integration.

\begin{lemma}(Strict Orthogonalization Error Bound)\label{lem:Strict Orthogonalization Error Bound}
Let $\Sigma$ be parameterized by $F(y)=u(y)e_1+y$ over the vertical hyperplane $\Pi\subset\mathbb{R}^{n+1}$. Let $\Omega\subset\Pi$ be the exterior domain $\Pi\setminus B_R(0)$ and let $\mu$ be the outward unit normal to $\partial B_R\subset\Pi$.

Define the pushforward vector $X=dF(\mu)=\mu+(D_\mu u)e_1$. Let $\eta\in T\Sigma$ be the unique unit conormal vector field orthogonal to $\partial\Sigma_R$ such that $\langle\eta,X\rangle>0$. Assume $\lim_{|y|\to\infty}|Du(y)|=0$. Then the error vector $v_{\mathrm{err}}:=\eta-X$ satisfies the asymptotic bound $v_{\mathrm{err}}\in\mathcal{B}_\infty(|Du|^2)$.
\end{lemma}

\begin{proof}
Let $\{\tau_1,\dots,\tau_{n-1}\}$ be an orthonormal basis for $T_y(\partial B_R)$. The tangent space $T_{F(y)}(\partial\Sigma_R)$ is spanned by the pushforward vectors $T_\alpha=\tau_\alpha+(D_{\tau_\alpha}u)e_1$. The induced metric tensor is $\gamma_{\alpha\beta}=\langle T_\alpha,T_\beta\rangle=\delta_{\alpha\beta}+(D_{\tau_\alpha}u)(D_{\tau_\beta}u)$. In matrix form this can be written as $\gamma=I+aa^T$, where $a$ is the vector with components $a_\alpha=D_{\tau_\alpha}u$.

The eigenvalues of $\gamma$ are $1$ with multiplicity $n-2$ and $1+|D_Tu|^2$ with multiplicity $1$. Hence $\gamma$ is positive definite and the inverse matrix $\gamma^{-1}$ is also positive definite with maximal eigenvalue equal to $1$ by Lemma~\ref{lem:Global Uniform Ellipticity}.

The vector $X$ decomposes orthogonally as $X=W+P$, where $P\in T_{F(y)}(\partial\Sigma_R)$ and $W\perp T_{F(y)}(\partial\Sigma_R)$. The projection $P$ onto $T_{F(y)}(\partial\Sigma_R)$ is given by
$
P=\sum_{\alpha,\beta}(\gamma^{-1})_{\alpha\beta}\langle X,T_\beta\rangle T_\alpha.
$
Let $v$ denote the vector with components $v_\beta=\langle X,T_\beta\rangle$. Because $\mu\perp\tau_\beta$ in $\Pi$, we compute
$
v_\beta=\langle \mu+(D_\mu u)e_1,\tau_\beta+(D_{\tau_\beta}u)e_1\rangle=(D_\mu u)(D_{\tau_\beta}u).
$

The squared Euclidean norm of the projection satisfies $||P||^2=v^T\gamma^{-1}v$. Since the maximal eigenvalue of $\gamma^{-1}$ does not exceed $1$, it follows that $||P||^2\le||v||^2$. Evaluating $||v||^2$ gives
\begin{equation} 
||v||^2=\sum_{\beta=1}^{n-1}(D_\mu u)^2(D_{\tau_\beta}u)^2=(D_\mu u)^2|D_Tu|^2\le|Du|^2|Du|^2=|Du|^4.
\end{equation}
Taking square roots yields $||P||\le|Du|^2$.

By orthogonality of the decomposition $X=W+P$, one has $||W||^2=||X||^2-||P||^2$. Since $||X||^2=1+(D_\mu u)^2$, define $\omega=(D_\mu u)^2-||P||^2$ so that $||W||=\sqrt{1+\omega}$. Using the bound on $P$, the triangle inequality gives $|\omega|\le|Du|^2+|Du|^4$.

Because $\lim_{|y|\to\infty}|Du|=0$, there exists a radius $R_0>0$ such that $|Du(y)|^2\le\frac14$ whenever $|y|>R_0$. Hence $|\omega|\le\frac14+\frac1{16}=\frac5{16}<\frac12$. For $\omega\in(-\tfrac12,\tfrac12)$, Taylor's theorem with Lagrange remainder yields $(1+\omega)^{-1/2}=1+\epsilon$, where $\epsilon=-\frac12\omega+\frac38(1+\xi)^{-5/2}\omega^2$ for some $\xi$ between $0$ and $\omega$.

Since $\xi>-1/2$, one has $(1+\xi)^{-5/2}<(1/2)^{-5/2}=4\sqrt2<6$. Therefore
\begin{equation} 
|\epsilon|\le\frac12|\omega|+\frac94|\omega|^2\le\frac12(|Du|^2+|Du|^4)+\frac94(|Du|^2+|Du|^4)^2.
\end{equation}
Because $|Du|^2\le1/4$, higher-order terms are absorbed and there exists a constant $C_\epsilon>0$ such that $|\epsilon|\le C_\epsilon|Du|^2$.

The unit conormal is $\eta=W/||W||$. Substituting $W=X-P$ and $||W||^{-1}=1+\epsilon$ gives $\eta=(X-P)(1+\epsilon)=X+\epsilon X-P-\epsilon P$. The error vector therefore satisfies $v_{\mathrm{err}}=\eta-X=\epsilon X-(1+\epsilon)P$.

Using the triangle inequality together with the bounds $||X||\le\sqrt{1+\tfrac14}=\tfrac{\sqrt5}{2}$, $||P||\le|Du|^2$, and $|\epsilon|\le C_\epsilon|Du|^2$, we obtain
\begin{equation} 
||v_{\mathrm{err}}||\le|\epsilon|\,||X||+|1+\epsilon|\,||P||\le\left(C_\epsilon\frac{\sqrt5}{2}\right)|Du|^2+(1+C_\epsilon/4)|Du|^2.
\end{equation}
Factoring $|Du|^2$ shows $||v_{\mathrm{err}}||\le C_{\mathrm{final}}|Du|^2$ for all $|y|>R_0$. Hence $v_{\mathrm{err}}\in\mathcal{B}_\infty(|Du|^2)$.
\end{proof}

\begin{remark}
The discrepancy between the intrinsic conormal field $\eta$ and the pushforward of the extrinsic radial normal $X$ is structurally governed by the gradient. The error vector is proven to strictly reside in the asymptotic class $\mathcal{B}_\infty(|Du|^2)$, thereby ensuring that asymptotic orthogonality is preserved up to a rigorously controlled quadratic remainder.
\end{remark}

\subsection{Expansion of the Flux Integrand}

Having bounded the conormal error, we proceed to compute the explicit algebraic expansion of the flux density scalar product. By taking the inner product of the conormal vector field with the translation vector, we bifurcate the integrand into a principal geometric quotient (specifically, $\frac{y_n}{R}$) and a strict scalar remainder function. This algebraic separation is structurally required to subsequently exploit the orthogonal symmetries of the ambient measure space.

\begin{lemma}[Asymptotic Expansion of the Flux Term]\label{lem:Asymptotic Expansion of the Flux Term}
Let the ambient Euclidean space be orthogonally decomposed as $\mathbb{R}^{n+1}=\mathrm{span}(e_1)\oplus\Pi$, where $\Pi$ is a vertical hyperplane with basis $\{\hat e_1,\dots,\hat e_n\}$. Let the translation vector be $e_{n+1}=\hat e_n$. Let the hypersurface $\Sigma$ be parameterized by the embedding $F(y)=u(y)e_1+y$ for $y\in\Pi$. Let $\partial\Sigma_R=F(\partial B_R)$ denote the boundary of the surface over the ball $B_R\subset\Pi$. Let $\eta\in T\Sigma$ be the outward unit conormal to $\partial\Sigma_R$.

Then the flux integrand satisfies
\begin{equation}
\langle\eta,e_{n+1}^\top\rangle=\frac{y_n}{R}+E(y),
\end{equation}
where $\frac{y_n}{R}$ is the normal component to the sphere $\partial B_R$, and the remainder function $E(y)$ satisfies strict asymptotic boundedness $E\in\mathcal B_\infty(|Du|^2)$.
\end{lemma}

\begin{proof}
The vector $e_{n+1}\in\mathbb R^{n+1}$ admits the orthogonal decomposition $e_{n+1}=e_{n+1}^\top+e_{n+1}^\perp$ where $e_{n+1}^\top\in T\Sigma$ and $e_{n+1}^\perp\in N\Sigma$. Since $\eta\in T\Sigma$, one has $\langle\eta,e_{n+1}^\perp\rangle=0$. By bilinearity of the Euclidean inner product,
\begin{equation}
\langle\eta,e_{n+1}^\top\rangle
=\langle\eta,e_{n+1}-e_{n+1}^\perp\rangle
=\langle\eta,e_{n+1}\rangle-\langle\eta,e_{n+1}^\perp\rangle
=\langle\eta,e_{n+1}\rangle
=\langle\eta,\hat e_n\rangle .
\end{equation}

Let $\mu(y)=y/R$ be the outward unit normal to the sphere $\partial B_R$ in the hyperplane $\Pi$. By Lemma~\ref{lem:Strict Orthogonalization Error Bound}, the conormal vector admits the decomposition $\eta=\mu+(D_\mu u)e_1+v_{err}$ where $v_{err}\in\mathcal B_\infty(|Du|^2)$.

Substituting this decomposition into the inner product gives
\begin{equation}
\langle\eta,\hat e_n\rangle
=\langle\mu+(D_\mu u)e_1+v_{err},\hat e_n\rangle
=\langle\mu,\hat e_n\rangle+\langle(D_\mu u)e_1,\hat e_n\rangle+\langle v_{err},\hat e_n\rangle .
\end{equation}

The vector $\mu$ has the representation $\mu=\sum_{k=1}^n (y_k/R)\hat e_k$. Since the basis is orthonormal, $\langle\mu,\hat e_n\rangle=y_n/R$. The vector $e_1$ spans the orthogonal complement of $\Pi$, while $\hat e_n\in\Pi$, hence $\langle e_1,\hat e_n\rangle=0$, which implies $\langle(D_\mu u)e_1,\hat e_n\rangle=0$. Therefore
\begin{equation}
\langle\eta,e_{n+1}^\top\rangle=\frac{y_n}{R}+\langle v_{err},\hat e_n\rangle .
\end{equation}

Define the scalar remainder function $E(y)=\langle v_{err}(y),\hat e_n\rangle$. By the Cauchy--Schwarz inequality,
\begin{equation}
|E(y)|=|\langle v_{err}(y),\hat e_n\rangle|\le \|v_{err}(y)\|\cdot\|\hat e_n\|=\|v_{err}(y)\|.
\end{equation}

Since $v_{err}\in\mathcal B_\infty(|Du|^2)$, there exist a constant $C>0$ and a compact set $K$ such that $\|v_{err}(y)\|\le C|Du(y)|^2$ for all $y\in\Pi\setminus K$. Consequently $|E(y)|\le C|Du(y)|^2$, which shows $E\in\mathcal B_\infty(|Du|^2)$. This proves the claim.
\end{proof}

\subsection{Measure Pullback and Pullback Flux}

To perform rigorous measure-theoretic integration, the flux form defined on the immersed manifold boundary must be formally pulled back to the standard Euclidean sphere. We establish the exact Jacobian distortion between the intrinsic Hausdorff measure and the standard Euclidean spherical measure. By consolidating the metric distortion and the integrand remainder, we define a unified pullback flux density, effectively reducing the geometric integral to a parameterized form on the standard sphere governed by strict asymptotic error bounds.

\begin{lemma}[Measure Pullback]\label{lem:Measure Pullback}
Let $F_k : \Pi_k \setminus \Omega_k \to \mathbb{R}^{n+1}$ be the immersion defined by $F_k(y)=u_k(y)e_1+y$. Let $\mathcal{S}_R \subset \Pi_k$ be the Euclidean $(n-1)$--sphere of radius $R$. Let $d\mathcal{H}^{n-1}_\Sigma$ denote the intrinsic Hausdorff measure on the embedded boundary $\partial E_{k,R}=F_k(\mathcal{S}_R)$ and let $d\sigma$ denote the standard Euclidean measure on $\mathcal{S}_R$. Then
\begin{equation} 
d\mathcal{H}^{n-1}_\Sigma=\sqrt{1+\|D_\tau u_k\|^2}\,d\sigma
\end{equation}
where $D_\tau u_k$ denotes the tangential projection of the ambient gradient $Du_k$ onto $T\mathcal{S}_R$.  

If $\lim_{|y|\to\infty}\|Du_k(y)\|=0$, there exists a function $J_{err}\in\mathcal{B}_\infty(|Du_k|^2)$ such that
\begin{equation} 
d\mathcal{H}^{n-1}_\Sigma=(1+J_{err}(y))\,d\sigma .
\end{equation}
\end{lemma}

\begin{proof}
Let $y\in\mathcal{S}_R$. The tangent space $T_y\Pi_k\cong\mathbb{R}^n$ decomposes orthogonally as $T_y\mathcal{S}_R\oplus N_y\mathcal{S}_R$. Let $\{\tau_1,\dots,\tau_{n-1}\}$ be a local orthonormal frame for $T_y\mathcal{S}_R$.

Consider the restricted immersion $\tilde F_k=F_k|_{\mathcal{S}_R}$. Its differential acts on the basis vectors as $d\tilde F_k(\tau_\alpha)=\tau_\alpha+(D_{\tau_\alpha}u_k)e_1$ for $\alpha=1,\dots,n-1$, where $D_{\tau_\alpha}u_k$ denotes the directional derivative along $\tau_\alpha$.

The induced first fundamental form on $\partial E_{k,R}$ has components $g_{\alpha\beta}=\langle d\tilde F_k(\tau_\alpha),d\tilde F_k(\tau_\beta)\rangle$. Substituting the differential gives $g_{\alpha\beta}=\langle \tau_\alpha+(D_{\tau_\alpha}u_k)e_1,\tau_\beta+(D_{\tau_\beta}u_k)e_1\rangle$. Because the vectors $\tau_\alpha$ are orthonormal in $\Pi_k$ and orthogonal to $e_1$, the mixed terms vanish and we obtain $g_{\alpha\beta}=\delta_{\alpha\beta}+(D_{\tau_\alpha}u_k)(D_{\tau_\beta}u_k)$.

In matrix form the metric tensor is $g=I_{n-1}+vv^T$, where $v=(D_{\tau_1}u_k,\dots,D_{\tau_{n-1}}u_k)^T$. The determinant of this rank--one perturbation satisfies $\det(I_{n-1}+vv^T)=1+v^Tv$. Therefore $\det(g)=1+\|D_\tau u_k\|^2$.

The area formula for immersions gives $d\mathcal{H}^{n-1}_\Sigma=\sqrt{\det(g)}\,d\sigma$. Substituting the determinant yields $d\mathcal{H}^{n-1}_\Sigma=\sqrt{1+\|D_\tau u_k\|^2}\,d\sigma$.

Define $J_{err}(y)=\sqrt{1+\|D_\tau u_k(y)\|^2}-1$. Assume $\lim_{|y|\to\infty}\|Du_k(y)\|=0$. Then there exists $R_{analytic}>0$ such that for all $|y|>R_{analytic}$ one has $\|Du_k(y)\|^2\le\frac12$. Because the tangential gradient is a projection, $\|D_\tau u_k\|^2\le\|Du_k\|^2\le\frac12$.

Let $x=\|D_\tau u_k\|^2$. Since $x\in[0,\tfrac12]$, the function $\sqrt{1+x}$ admits the expansion $\sqrt{1+x}=1+\frac12x-\frac{1}{8(1+\xi)^{3/2}}x^2$ for some $\xi\in(0,x)$. In particular $\sqrt{1+x}-1\le\frac12x$ for all $x\ge0$. Hence $0\le J_{err}(y)\le\frac12\|D_\tau u_k\|^2\le\frac12\|Du_k\|^2$.

Thus for all $|y|>R_{analytic}$ one has $|J_{err}(y)|\le C'\|Du_k\|^2$ with $C'=\frac12$. This shows $J_{err}\in\mathcal{B}_\infty(|Du_k|^2)$ and completes the proof.
\end{proof}

\begin{proposition}[Boundary Flux Integrand]\label{prop:Boundary Flux Integrand}
Let $F_k : \Pi_k \setminus \Omega_k \to \mathbb{R}^{n+1}$ be the immersion parameterizing the translator end $E_k$. Let $\mathcal{S}_R = \partial B_R(0) \subset \Pi_k$. Let $\omega_{flux} = \langle \eta_k, e_{n+1}^\top \rangle d\mathcal{H}^{n-1}_\Sigma$ denote the boundary flux density on $\partial E_{k,R}$. Under the pullback $F_k^*$ the form satisfies the asymptotic expansion on $\mathcal{S}_R$
\begin{equation} 
F_k^* \omega_{flux} = \left(\frac{y_n}{R} + \tilde{E}_k(y)\right) d\sigma
\end{equation}
where the remainder $\tilde{E}_k$ satisfies $\tilde{E}_k \in \mathcal{B}_\infty(|Du_k|^2)$.
\end{proposition}

\begin{proof}
By the definition of pullback for densities, the pullback form $F_k^* \omega_{flux}$ is obtained by evaluating the scalar factor of $\omega_{flux}$ along the immersion and multiplying it by the pullback of the measure.

From the pointwise computation of the boundary normal (Lemma~\ref{lem:Asymptotic Expansion of the Flux Term}) we have
\begin{equation} 
\langle \eta_k, e_{n+1}^\top \rangle \circ F_k(y) = \frac{y_n}{R} + E_k(y)
\end{equation}
where $E_k \in \mathcal{B}_\infty(|Du_k|^2)$.  

From the pullback of the boundary measure (Lemma~\ref{lem:Measure Pullback}) one obtains
\begin{equation} 
F_k^* d\mathcal{H}^{n-1}_\Sigma = (1 + J_{err}(y)) d\sigma
\end{equation}
with $J_{err} \in \mathcal{B}_\infty(|Du_k|^2)$.

Substituting these identities into the pullback of the flux form gives
\begin{equation} 
F_k^* \omega_{flux} = \left(\frac{y_n}{R} + E_k(y)\right)(1 + J_{err}(y)) d\sigma .
\end{equation}

Expanding the scalar product yields
\begin{equation} 
F_k^* \omega_{flux}
=
\left(
\frac{y_n}{R}
+
\frac{y_n}{R} J_{err}(y)
+
E_k(y)
+
E_k(y) J_{err}(y)
\right) d\sigma .
\end{equation}

Define the unified remainder
\begin{equation} 
\tilde{E}_k(y)
=
\frac{y_n}{R} J_{err}(y)
+
E_k(y)
+
E_k(y) J_{err}(y) .
\end{equation}

We verify that each term belongs to $\mathcal{B}_\infty(|Du_k|^2)$.  
On $\mathcal{S}_R$ the coordinate $y_n$ satisfies $|y_n| \le R$, hence the factor $\frac{y_n}{R}$ is uniformly bounded. Since $J_{err} \in \mathcal{B}_\infty(|Du_k|^2)$, bounded multiplication implies $\frac{y_n}{R} J_{err}(y) \in \mathcal{B}_\infty(|Du_k|^2)$.

The term $E_k(y)$ already satisfies $E_k \in \mathcal{B}_\infty(|Du_k|^2)$.

For the product term $E_k(y) J_{err}(y)$ there exist constants $C_1$ and $C_2$ such that $|E_k(y)| \le C_1 |Du_k|^2$ and $|J_{err}(y)| \le C_2 |Du_k|^2$ outside a compact set. Hence
\begin{equation} 
|E_k(y) J_{err}(y)| \le (C_1 C_2) |Du_k|^4 .
\end{equation}
Thus this term lies in $\mathcal{B}_\infty(|Du_k|^4)$. Because $\lim_{|y|\to\infty} \|Du_k(y)\| = 0$, the inequality $|Du_k|^4 \le |Du_k|^2$ holds outside a sufficiently large compact set, and therefore the product term also belongs to $\mathcal{B}_\infty(|Du_k|^2)$.

Since the class $\mathcal{B}_\infty(|Du_k|^2)$ is closed under addition, the sum defining $\tilde{E}_k$ satisfies $\tilde{E}_k \in \mathcal{B}_\infty(|Du_k|^2)$.

Substituting this definition into the expansion gives
\begin{equation} 
F_k^* \omega_{flux} = \left(\frac{y_n}{R} + \tilde{E}_k(y)\right) d\sigma .
\end{equation}
\end{proof}

\begin{remark}
Through the algebraic pullback of the boundary density onto the Euclidean sphere $\mathcal{S}_R$, the flux form is canonically decomposed. The principal anisotropic term $\frac{y_n}{R}$ is strategically isolated, while all higher-order geometric non-linearities—arising from the metric and the conormal error—are strictly absorbed into a unified scalar remainder $\tilde{E}_k$ that inherits the quadratic decay class of the gradient.
\end{remark}

\section{Main Result}

This final chapter synthesizes the previously constructed analytic machinery and geometric identities to formulate the ultimate structural classification of the translator. By establishing the vanishing of the flux limit, we strictly deduce the global triviality of the mean curvature, culminating in the proof of geometric rigidity.

\subsection{Liouville Rigidity for Minimal Graphs}

Before executing the global limit, we must isolate the analytic rigidity of the trivial mean curvature state. We prove that if the mean curvature vanishes identically and the parameterization obeys asymptotic gradient triviality, the resulting minimal surface equation structurally reduces to a globally uniformly elliptic operator. By applying the weak maximum principle to the differentiated minimal surface equation, we establish a Liouville-type rigidity theorem, demonstrating that such a graph is strictly constant and hence forms a vertical hyperplane.

\begin{lemma}(Global Uniform Ellipticity)\label{lem:Global Uniform Ellipticity}
Let $\tilde{u} \in C^2(\mathbb{R}^{n-1}, \mathbb{R})$ be a scalar function satisfying the strict asymptotic triviality condition $\lim_{|y'|\to\infty} ||D\tilde{u}(y')|| = 0$. Let the linear operator coefficients be defined by:
\begin{equation}
a_{ij}(D\tilde{u}) = \frac{1}{\sqrt{1+||D\tilde{u}||^2}} \left( \delta_{ij} - \frac{D_i \tilde{u} D_j \tilde{u}}{1+||D\tilde{u}||^2} \right)
\end{equation}
Then the coefficient matrix is globally uniformly elliptic. Specifically, there exist universal constants $0 < \lambda_{min} \le \lambda_{max} < \infty$ such that for all $y' \in \mathbb{R}^{n-1}$ and all $\xi \in \mathbb{R}^{n-1}$:
\begin{equation}
\lambda_{min} ||\xi||^2 \le \sum_{i,j=1}^{n-1} a_{ij}(D\tilde{u}(y')) \xi_i \xi_j \le \lambda_{max} ||\xi||^2
\end{equation}
\end{lemma}

\begin{proof}
We evaluate the quadratic form $Q(y', \xi)$ associated with the matrix $a_{ij}(D\tilde{u})$ at an arbitrary point $y' \in \mathbb{R}^{n-1}$ and an arbitrary test vector $\xi \in \mathbb{R}^{n-1}$. Let $p = D\tilde{u}(y')$.
\begin{equation}
Q(y', \xi) = \sum_{i,j=1}^{n-1} a_{ij}(p) \xi_i \xi_j = \frac{1}{\sqrt{1+||p||^2}} \left( ||\xi||^2 - \frac{\langle p, \xi \rangle^2}{1+||p||^2} \right)
\end{equation}

By the fundamental positivity of real squares, $\frac{\langle p, \xi \rangle^2}{1+||p||^2} \ge 0$. Therefore, the term inside the parenthesis is strictly bounded above by $||\xi||^2$.

Furthermore, because $||p||^2 \ge 0$, the scaling factor satisfies $(1+||p||^2)^{-1/2} \le 1$.

Consequently:
\begin{equation}
Q(y', \xi) \le 1 \cdot ||\xi||^2
\end{equation}
This establishes the strict upper ellipticity constant $\lambda_{max} = 1$, valid globally across $\mathbb{R}^{n-1}$.

To bound the quadratic form from below, we apply the Cauchy-Schwarz inequality to the standard Euclidean inner product: $\langle p, \xi \rangle^2 \le ||p||^2 ||\xi||^2$.

Substituting this into the subtractive term of the quadratic form yields:
\begin{equation}
Q(y', \xi) \ge \frac{1}{\sqrt{1+||p||^2}} \left( ||\xi||^2 - \frac{||p||^2 ||\xi||^2}{1+||p||^2} \right)
\end{equation}
Algebraically factoring $||\xi||^2$ from the parenthesis produces:
\begin{equation}
Q(y', \xi) \ge \frac{||\xi||^2}{\sqrt{1+||p||^2}} \left( 1 - \frac{||p||^2}{1+||p||^2} \right) = \frac{||\xi||^2}{\sqrt{1+||p||^2}} \left( \frac{1}{1+||p||^2} \right)
\end{equation}
Combining the exponents, the strict algebraic lower bound for the minimal eigenvalue at point $y'$ is exactly:
\begin{equation}
\lambda(y') \ge \frac{1}{(1+||D\tilde{u}(y')||^2)^{3/2}}
\end{equation}

Because $\tilde{u} \in C^2(\mathbb{R}^{n-1})$, its gradient vector field $D\tilde{u}$ is continuous on $\mathbb{R}^{n-1}$. By hypothesis, $\lim_{|y'|\to\infty} ||D\tilde{u}(y')|| = 0$.Hence the global supremum of the gradient magnitude is strictly finite. Let $M > 0$ be this finite maximum:
\begin{equation}
\sup_{y' \in \mathbb{R}^{n-1}} ||D\tilde{u}(y')|| = M < \infty
\end{equation}
Because the function $f(x) = (1+x^2)^{-3/2}$ is strictly monotonically decreasing for $x \ge 0$, replacing $||D\tilde{u}(y')||$ with its global maximum $M$ strictly minimizes the algebraic lower bound over the entire domain:
\begin{equation}
\frac{1}{(1+||D\tilde{u}(y')||^2)^{3/2}} \ge \frac{1}{(1+M^2)^{3/2}}
\end{equation}
We define the universal constant $\lambda_{min} = (1+M^2)^{-3/2}$. Because $M$ is finite, $\lambda_{min}$ is strictly positive ($\lambda_{min} > 0$).

Thus, for all $y' \in \mathbb{R}^{n-1}$, we possess the  global bound:
\begin{equation}
Q(y', \xi) \ge \lambda_{min} ||\xi||^2
\end{equation}
This completes the proof.
\end{proof}

\begin{theorem}[Liouville Rigidity for Asymptotic Minimal Graphs]\label{thm:Liouville Rigidity for Asymptotic Minimal Graphs}
Let $\Omega \subset \mathbb{R}^n$ be a bounded open set and let $u:\mathbb{R}^n \setminus \Omega \to \mathbb{R}$ be a smooth function whose graph $\Sigma=\mathrm{graph}(u)$ is a hypersurface satisfying $H\equiv0$. Assume that the projection of the vertical translation vector onto the normal bundle vanishes, $e_{n+1}^{\perp}=0$, and that $\lim_{|y|\to\infty}\|Du(y)\|=0$. Then $u$ is constant and $\Sigma$ is a vertical hyperplane.
\end{theorem}

\begin{proof}
Write the ambient coordinates as $(x_1,y_1,\dots,y_n)$ and denote by $\hat e_1,\dots,\hat e_n$ the coordinate vectors in the $y$–directions, with $e_{n+1}=\hat e_n$. The graph $x_1=u(y)$ has unit normal
\begin{equation}
\nu=\frac{e_1-\sum_{i=1}^n D_i u\,\hat e_i}{\sqrt{1+\|Du\|^2}} .
\end{equation}
The condition $e_{n+1}^{\perp}=0$ is equivalent to $\langle e_{n+1},\nu\rangle=0$. Substituting the expression for $\nu$ gives
\begin{equation}
\langle \hat e_n,\nu\rangle=-\frac{D_n u}{\sqrt{1+\|Du\|^2}}=0 .
\end{equation}
Hence $D_n u\equiv0$ on $\mathbb{R}^n\setminus\Omega$, so $u$ is independent of $y_n$.

Because $\Omega$ is bounded, there exists $M>0$ with $\Omega\subset B_M(0)$. For any $y'=(y_1,\dots,y_{n-1})\in\mathbb{R}^{n-1}$ and any $|y_n|>M$, the point $(y',y_n)$ lies outside $\Omega$. Since $u$ does not depend on $y_n$, define $\tilde u:\mathbb{R}^{n-1}\to\mathbb{R}$ by $\tilde u(y')=u(y',y_n)$ for any $|y_n|>M$. This defines a smooth entire function on $\mathbb{R}^{n-1}$. The decay assumption implies $\lim_{|y'|\to\infty}\|D\tilde u(y')\|=0$.

Because $H\equiv0$, the function $\tilde u$ satisfies the minimal surface equation
\begin{equation}
\sum_{i=1}^{n-1} D_i\!\left(\frac{D_i\tilde u}{\sqrt{1+\|D\tilde u\|^2}}\right)=0 .
\end{equation}
Let $v_k=D_k\tilde u$ for $1\le k\le n-1$. Differentiating the equation with respect to $y_k$ yields
\begin{equation}
\sum_{i,j=1}^{n-1} D_i\!\left(a_{ij}(D\tilde u)\,D_j v_k\right)=0 ,
\end{equation}
where $a_{ij}(p)=(1+\|p\|^2)^{-1/2}\bigl(\delta_{ij}-(1+\|p\|^2)^{-1}p_i p_j\bigr)$. The matrix $(a_{ij})$ is positive definite, so each $v_k$ satisfies a uniformly elliptic linear equation in divergence form with no zeroth-order term.

Let $\varepsilon>0$. From $\lim_{|y'|\to\infty}\|D\tilde u(y')\|=0$, there exists $R>0$ such that $|v_k(y')|\le\varepsilon$ whenever $|y'|>R$. On the compact ball $\overline{B_R}$ the function $v_k$ is continuous and satisfies the elliptic equation in the interior. By the weak maximum principle, the maximum and minimum of $v_k$ on $\overline{B_R}$ occur on $\partial B_R$. Since $|v_k|\le\varepsilon$ on $\partial B_R$, it follows that $|v_k(y')|\le\varepsilon$ for all $y'\in\overline{B_R}$. Therefore $|v_k(y')|\le\varepsilon$ on $\mathbb{R}^{n-1}$.

Because $\varepsilon$ was arbitrary, $v_k\equiv0$ for every $k$, so $D\tilde u\equiv0$. Hence $\tilde u$ is constant. Consequently $u$ is constant on $\mathbb{R}^n\setminus\Omega$, and its graph is a vertical hyperplane. If the asymptotic plane is normalized to $x_1=0$, the constant is zero and $u\equiv0$.
\end{proof}

\begin{remark}
The intersection of the minimal surface condition $H \equiv 0$ and the orthogonality of the translation vector $e_{n+1}^\perp = 0$ strictly suppresses the dependence on the $y_n$ coordinate. This reduction to a lower-dimensional minimal surface equation triggers a classical Liouville-type rigidity, which strictly forces any solution possessing an asymptotically trivial gradient to be a global constant, manifesting geometrically as a flat vertical hyperplane.
\end{remark}

\subsection{Refined Remainder Estimates and Spherical Integrals}

To evaluate the limit of the pullback flux, we must govern the unified remainder term over the expanding spherical domain. We apply the bifurcated, pointwise gradient decay bounds established via the spatial partitions (the Wake and Exponential sets) to strictly bound the remainder function. We formalize the integration over these measurable spherical subsets, yielding strict, radius-dependent algebraic decay rates for the subsequent flux analysis.

\begin{lemma}[Bifurcated Remainder Estimate]\label{lem:Bifurcated Remainder Estimate}
Let $\Pi_{k} \cong \mathbb{R}^{n}$ be a Euclidean hyperplane. Let $\tilde{E}_{k}: \Pi_{k} \to \mathbb{R}$ satisfy strict asymptotic boundedness relative to the squared gradient, denoted $\tilde{E}_{k} \in \mathcal{B}_{\infty}(|Du_{k}|^{2})$.

Assume the gradient satisfies the bifurcated pointwise decay properties as in Theorem~\ref{thm:priori estimate and gradient bifurcation}:

\begin{equation}
|Du_{k}(y)|^{2} \le C_{\mathcal{W}}|y|^{-n} \quad \text{for } y \in \mathcal{W}_{C} \setminus K_{\mathcal{W}},
\end{equation}
\begin{equation}
|Du_{k}(y)|^{2} \le C_{\mathcal{E}}e^{-2\mu |y|} \quad \text{for } y \in \mathcal{E}_{C} \setminus K_{\mathcal{E}},
\end{equation}
where $C_{\mathcal{W}}, C_{\mathcal{E}}, \mu$ are strictly positive real constants, and $K_{\mathcal{W}}, K_{\mathcal{E}}$ are compact subsets of $\Pi_{k}$.

Then there exist positive real constants $M_{\mathcal{W}}$ and $M_{\mathcal{E}}$, and a global critical radius $R_{crit} > 0$, such that for all $y \in \Pi_{k}$ satisfying $|y| \ge R_{crit}$:

\begin{enumerate}
    \item If $y \in \mathcal{W}_{C}$, then $|\tilde{E}_{k}(y)| \le M_{\mathcal{W}}|y|^{-n}$.
    \item If $y \in \mathcal{E}_{C}$, then $|\tilde{E}_{k}(y)| \le M_{\mathcal{E}}e^{-2\mu |y|}$.
\end{enumerate}
\end{lemma}

\begin{proof}
By Definition,, the hypothesis $\tilde{E}_{k} \in \mathcal{B}_{\infty}(|Du_{k}|^{2})$ implies the existence of a strictly positive real constant $C_{3} > 0$ and a compact set $K_{3} \subset \Pi_{k}$ such that for all $y \in \Pi_{k} \setminus K_{3}$, the strict algebraic inequality
\begin{equation}
|\tilde{E}_{k}(y)| \le C_{3}|Du_{k}(y)|^{2}
\end{equation}
holds.

We construct the universal exclusion domain. Let $K_{universal} := K_{\mathcal{W}} \cup K_{\mathcal{E}} \cup K_{3}$, then $K_{universal}$ is compact, and there exists a strictly positive radius $R_{crit} > 0$ such that $K_{universal} \subset B_{R_{crit}}(0)$.

Consequently, for any $y \in \Pi_{k}$, the metric condition $|y| \ge R_{crit}$  implies $y \notin B_{R_{crit}}(0)$, which strictly dictates $y \notin K_{universal}$. By the definition of the union, $y \notin K_{universal}$ implies $y \notin K_{\mathcal{W}}$, $y \notin K_{\mathcal{E}}$, and $y \notin K_{3}$.

Let $y \in \Pi_{k}$ such that $|y| \ge R_{crit}$. Because $y \notin K_{3}$, the asymptotic bound is  active:
\begin{equation}
|\tilde{E}_{k}(y)| \le C_{3}|Du_{k}(y)|^{2}.
\end{equation}

We evaluate this globally active bound over the strictly defined topological partition:

\textbf{Case 1: $y \in \mathcal{W}_{C}$}.  
Because $y \notin K_{\mathcal{W}}$, the polynomial decay hypothesis is active: 
\begin{equation}
|Du_{k}(y)|^{2} \le C_{\mathcal{W}}|y|^{-n}.
\end{equation}
Substituting this upper bound into the asymptotic bound yields:
\begin{equation}
|\tilde{E}_{k}(y)| \le C_{3}(C_{\mathcal{W}}|y|^{-n}).
\end{equation}
Defining the constant $M_{\mathcal{W}} := C_{3}C_{\mathcal{W}}$, we obtain the  strict inequality
\begin{equation}
|\tilde{E}_{k}(y)| \le M_{\mathcal{W}}|y|^{-n}.
\end{equation}

\textbf{Case 2: $y \in \mathcal{E}_{C}$}.  
Because $y \notin K_{\mathcal{E}}$, the exponential decay hypothesis is active:
\begin{equation}
|Du_{k}(y)|^{2} \le C_{\mathcal{E}}e^{-2\mu |y|}.
\end{equation}
Substituting this upper bound into the asymptotic bound yields:
\begin{equation}
|\tilde{E}_{k}(y)| \le C_{3}(C_{\mathcal{E}}e^{-2\mu |y|}).
\end{equation}
Defining the constant $M_{\mathcal{E}} := C_{3}C_{\mathcal{E}}$, we obtain the  strict inequality
\begin{equation}
|\tilde{E}_{k}(y)| \le M_{\mathcal{E}}e^{-2\mu |y|}.
\end{equation}

This exhaustive evaluation over the partition completes the proof.
\end{proof}

\begin{remark}
It is proper to note that the bifurcated pointwise decay of the gradient vector field (established in Theorem 3.4) is canonically inherited by the geometric remainder $\tilde{E}_k$. The strict substitution of the gradient bounds amplifies the spatial decay, rigorously confining the remainder to an $\mathcal{O}(|y|^{-n})$ profile in the wake set $\mathcal{W}_C$ and an $\mathcal{O}(e^{-2\mu|y|})$ profile in the exponential set $\mathcal{E}_C$.
\end{remark}

\begin{lemma}[Spherical Integral Bound]\label{lem:Spherical Integral Bound}
Let $\mathcal{S}_{R} = \partial B_{R}(0) \subset \mathbb{R}^{n}$ be the Euclidean sphere equipped with surface measure $d\sigma$. Let $A \subset \mathcal{S}_{R}$ be a Borel-measurable subset. For any non-negative measurable function $f: A \to \mathbb{R}_{\ge 0}$ satisfying a uniform upper bound $f(y) \le M$ for all $y \in A$, the integral  satisfies:
\begin{equation}
\int_{A} f(y) d\sigma \le M \omega_{n-1} R^{n-1}.
\end{equation}
\end{lemma}

\begin{proof}
By Definition, the local integral over $A$ is strictly equivalent to the global integral over the entire sphere $\mathcal{S}_{R}$ against the indicator function $\chi_{A}$:
\begin{equation}
\int_{A} f(y) d\sigma = \int_{\mathcal{S}_{R}} f(y) \chi_{A}(y) d\sigma.
\end{equation}
We establish a pointwise bound for the integrand $f(y) \chi_{A}(y)$ valid for all $y \in \mathcal{S}_{R}$.

By the definition of the indicator function, if $y \notin A$, then $\chi_{A}(y) = 0$. Since $f$ is non-negative, $f(y) \chi_{A}(y) = 0 \le M \cdot 0 = 0$.

If $y \in A$, then $\chi_{A}(y) = 1$. By the proposition's hypothesis, $f(y) \le M$. Multiplying by the strictly positive $\chi_{A}(y)$ yields $f(y) \chi_{A}(y) \le M \chi_{A}(y)$.

Thus, for all $y \in \mathcal{S}_{R}$, we possess the  strict inequality:
\begin{equation}
f(y) \chi_{A}(y) \le M \chi_{A}(y).
\end{equation}
Furthermore, by the definition of subsets, $A \subseteq \mathcal{S}_{R}$, which dictates that for all $y \in \mathcal{S}_{R}$, the indicator function satisfies $\chi_{A}(y) \le \chi_{\mathcal{S}_{R}}(y) = 1$. Since $M > 0$, multiplying this inequality by $M$ yields:
\begin{equation}
M \chi_{A}(y) \le M \cdot 1.
\end{equation}
By the transitivity of the real inequality, we conclude that for almost every $y \in \mathcal{S}_{R}$ (in fact, for all $y$):
\begin{equation}
f(y) \chi_{A}(y) \le M.
\end{equation}
By the strict monotonicity of the Lebesgue integral, integrating this uniformly bounded non-negative functiver the measure space $(\mathcal{S}_{R}, \sigma)$ preserves the inequality:
\begin{equation}
\int_{\mathcal{S}_{R}} f(y) \chi_{A}(y) d\sigma \le \int_{\mathcal{S}_{R}} M d\sigma.
\end{equation}
By the linearity of the integral, the constant $M$ factors out:
\begin{equation}
\int_{\mathcal{S}_{R}} M d\sigma = M \int_{\mathcal{S}_{R}} 1 d\sigma = M \sigma(\mathcal{S}_{R}).
\end{equation}
Finally, the total measure of the parameterized sphere strictly equals $\omega_{n-1} R^{n-1}$. Substituting this geometric magnitude into the inequality yields the final bound:
\begin{equation}
\int_{A} f(y) d\sigma \le M \omega_{n-1} R^{n-1}.
\end{equation}
\end{proof}

\subsection{Vanishing of the Boundary Flux}

In this section, we definitively evaluate the asymptotic limit of the boundary flux integral as the exhaustion radius diverges. Utilizing a formal orientation-reversing isometry, we demonstrate that the integral of the principal quotient strictly vanishes. Simultaneously, by invoking the bifurcated spherical bounds, we prove that the remainder integral is crushed in the limit by the dominant polynomial and exponential decay, rigorously establishing that the total boundary flux asymptotically approaches zero.

\begin{theorem}[Asymptotic Flux Vanishing]\label{thm:Asymptotic Flux Vanishing}
Let $n \ge 3$. Let $E_{k, R}$ be the cylindrical exhaustion of the translator end parameterized over the hyperplane $\Pi_{k} \cong \mathbb{R}^{n}$. Under the strict measurable partition defined by $\mathcal{W}_{C}$ and $\mathcal{E}_{C}$, the boundary flux integral  vanishes in the limit:
\begin{equation}
\lim_{R \to \infty} \int_{\partial E_{k, R}} \langle \eta_{k}, e_{n+1}^{\top} \rangle d\mathcal{H}_{\Sigma}^{n-1} = 0.
\end{equation}
\end{theorem}

\begin{proof}
By Proposition~\ref{prop:Boundary Flux Integrand}, the pullback of the boundary flux density onto the Euclidean sphere $\mathcal{S}_{R} \subset \Pi_{k}$ expands strictly to:
\begin{equation}
\int_{\partial E_{k, R}} \langle \eta_{k}, e_{n+1}^{\top} \rangle d\mathcal{H}_{\Sigma}^{n-1} = \int_{\mathcal{S}_{R}} \frac{y_{n}}{R} d\sigma + \int_{\mathcal{S}_{R}} \tilde{E}_{k}(y) d\sigma
\end{equation}
We define the principal term as 
\begin{equation}
I_{prin} := \int_{\mathcal{S}_{R}} \frac{y_{n}}{R} d\sigma
\end{equation}
and the remainder term as 
\begin{equation}
I_{rem}(R) := \int_{\mathcal{S}_{R}} \tilde{E}_{k}(y) d\sigma.
\end{equation}

We define the geometric reflection mapping $\Phi: \Pi_{k} \to \Pi_{k}$ by the strict coordinate assignment $\Phi(y_{1}, \dots, y_{n-1}, y_{n}) = (y_{1}, \dots, y_{n-1}, -y_{n})$. This mapping is an orthogonal involution, hence a global isometry, preserving both the sphere $\mathcal{S}_{R}$ and the Euclidean measure $d\sigma$.

Applying the change of variables formula under the pullback $\Phi^{*}$:
\begin{equation}
I_{prin} = \int_{\mathcal{S}_{R}} \Phi^{*}\left(\frac{y_{n}}{R} d\sigma\right) = \int_{\mathcal{S}_{R}} \frac{-y_{n}}{R} d\sigma = -I_{prin}
\end{equation}
The unique algebraic solution to $I_{prin} = -I_{prin}$ in $\mathbb{R}$ is $I_{prin} = 0$. This holds  for all $R > 0$.

We analyze $I_{rem}(R)$. By Lemma~\ref{lem:Spherical Measurable Partition}, for any $C > 0$, the sets $\mathcal{W}_{C, R} = \mathcal{W}_{C} \cap \mathcal{S}_{R}$ and $\mathcal{E}_{C, R} = \mathcal{E}_{C} \cap \mathcal{S}_{R}$ constitute a strict measurable partition of $\mathcal{S}_{R}$.

By the additivity of the Lebesgue integral over disjoint measurable domains, we rewrite the absolute value of the remainder:
\begin{equation}
|I_{rem}(R)| \le \int_{\mathcal{W}_{C, R}} |\tilde{E}_{k}(y)| d\sigma + \int_{\mathcal{E}_{C, R}} |\tilde{E}_{k}(y)| d\sigma
\end{equation}
Let $R > R_{crit}$ as established in Lemma~\ref{lem:Bifurcated Remainder Estimate},

On the measurable subset $\mathcal{W}_{C, R}$, every point satisfies $|y| = R$ and $|\tilde{E}_{k}(y)| \le M_{\mathcal{W}} R^{-n}$. Applying Lemma~\ref{lem:Spherical Integral Bound}, we establish the strict upper bound for the first integral:
\begin{equation}
\int_{\mathcal{W}_{C, R}} |\tilde{E}_{k}(y)| d\sigma \le (M_{\mathcal{W}} R^{-n}) \cdot (\omega_{n-1} R^{n-1}) = M_{\mathcal{W}} \omega_{n-1} R^{-1}.
\end{equation}
On the measurable subset $\mathcal{E}_{C, R}$, every point satisfies $|y| = R$ and $|\tilde{E}_{k}(y)| \le M_{\mathcal{E}} e^{-2\mu R}$. Applying Lemma~\ref{lem:Spherical Integral Bound} again, we establish the strict upper bound for the second integral:
\begin{equation}
\int_{\mathcal{E}_{C, R}} |\tilde{E}_{k}(y)| d\sigma \le (M_{\mathcal{E}} e^{-2\mu R}) \cdot (\omega_{n-1} R^{n-1}) = M_{\mathcal{E}} \omega_{n-1} R^{n-1} e^{-2\mu R}.
\end{equation}

Summing the localized bounds yields the global strict inequality for $R > R_{crit}$:
\begin{equation}
|I_{rem}(R)| \le M_{\mathcal{W}} \omega_{n-1} R^{-1} + M_{\mathcal{E}} \omega_{n-1} R^{n-1} e^{-2\mu R}.
\end{equation}
We evaluate the limits of both bounding terms as $R \to \infty$:
\begin{equation}
\lim_{R \to \infty} (M_{\mathcal{W}} \omega_{n-1} R^{-1}) = 0.
\end{equation}
Because the exponential decay strictly dominates polynomial growth for any degree $n-1$, 
\begin{equation}
\lim_{R \to \infty} (M_{\mathcal{E}} \omega_{n-1} R^{n-1} e^{-2\mu R}) = 0.
\end{equation}
By the Squeeze Theorem, 
\begin{equation}
\lim_{R \to \infty} |I_{rem}(R)| = 0,
\end{equation}
implying 
\begin{equation}
\lim_{R \to \infty} I_{rem}(R) = 0.
\end{equation}

Since 
\begin{equation}
\int_{\partial E_{k, R}} \langle \eta_{k}, e_{n+1}^{\top} \rangle d\mathcal{H}_{\Sigma}^{n-1} = I_{prin} + I_{rem}(R),
\end{equation}
and $I_{prin} = 0$ independently of $R$, taking the limit yields:
\begin{equation}
\lim_{R \to \infty} \int_{\partial E_{k, R}} \langle \eta_{k}, e_{n+1}^{\top} \rangle d\mathcal{H}_{\Sigma}^{n-1} = 0 + 0 = 0.
\end{equation}
\end{proof}

\begin{remark}
The vanishing of the boundary flux relies strictly on a precise geometric and measure-theoretic cancellation over the cylindrical exhaustion. The principal linear flux term evaluates exactly to zero due to the antipodal involution of the sphere, while the integral of the remainder terms is strictly suppressed in the limit by the bifurcated radial decay bounds acting against the spherical measure.
\end{remark}

\subsection{Global Rigidity of Asymptotic Translators}

We finalize the global classification theorem. By applying the zero limit of the boundary flux to the nested sequence of the universal flux identity, the monotone convergence theorem strictly forces the global $\mathcal{L}^2$ norm of the mean curvature to vanish. The identical vanishing of the mean curvature implies the immersion is a minimal hypersurface; bounded by our Liouville rigidity theorem and the principle of analytic continuation, the entire connected translator is globally forced into the trivial configuration of a vertical hyperplane.

\begin{theorem}[Rigidity of Asymptotic Translators]\label{thm:Rigidity of Asymptotic Translators}
Let $n \ge 3$. Let $\Sigma^n \hookrightarrow \mathbb{R}^{n+1}$ be a complete, properly immersed, connected translator with finite topology. Assume that every end of $\Sigma$ is $C^1$-asymptotic to a vertical hyperplane. Then $\Sigma$ is globally a vertical hyperplane.
\end{theorem}

\begin{proof}
Since $\Sigma$ has finite topology, there exists a compact core submanifold $K_{core} \subset \Sigma$ with smooth boundary such that $\Sigma \setminus K_{core}$ decomposes into finitely many disjoint ends $E_1,\dots,E_m$.

By the asymptotic hypothesis, each end $E_k$ admits a parameterization by a smooth immersion $F_k : \Pi_k \setminus \Omega_k \to \mathbb{R}^{n+1}$ of the form $F_k(y)=u_k(y)e_1+y$, where $\Pi_k$ is a vertical hyperplane, $\Omega_k \subset \Pi_k$ is bounded, and $\lim_{|y|\to\infty}\|Du_k(y)\|=0$.

For $R$ sufficiently large such that $\Omega_k \subset B_R(0) \subset \Pi_k$ for all $k$, define the truncated ends $E_{k,R}=F_k(B_R(0)\setminus \Omega_k)$. Define the global exhaustion
\begin{equation} 
\Sigma_{global,R}=K_{core}\cup\left(\bigcup_{k=1}^m E_{k,R}\right).
\end{equation}
This domain has boundary
\begin{equation} 
\partial\Sigma_{global,R}=\bigcup_{k=1}^m \partial E_{k,R}.
\end{equation}

Applying the universal flux identity (Theorem~\ref{thm:The Universal Flux Identity}) to $\Sigma_{global,R}$ gives
\begin{equation} 
\int_{\Sigma_{global,R}} |\vec{H}|^2\, d\mathcal{H}^n
=
\sum_{k=1}^m
\int_{\partial E_{k,R}}
\langle e_{n+1}^\top,\eta_k\rangle
\, d\mathcal{H}^{n-1}_\Sigma ,
\end{equation}
where $\eta_k$ denotes the outward unit conormal along $\partial E_{k,R}$.

For each end $k$, from Theorem~\ref{thm:Asymptotic Flux Vanishing} the asymptotic flux identity implies
\begin{equation} 
\lim_{R\to\infty}
\int_{\partial E_{k,R}}
\langle e_{n+1}^\top,\eta_k\rangle
\, d\mathcal{H}^{n-1}_\Sigma
=0.
\end{equation}
Since there are finitely many ends, the sum of these limits is also zero.

The domains $\Sigma_{global,R}$ form a nested sequence exhausting $\Sigma$. Because $|\vec{H}|^2 \ge 0$, the integrals $\int_{\Sigma_{global,R}} |\vec{H}|^2 d\mathcal{H}^n$ form a monotone increasing sequence. By the monotone convergence theorem,
\begin{equation} 
\lim_{R\to\infty}
\int_{\Sigma_{global,R}} |\vec{H}|^2\, d\mathcal{H}^n
=
\int_\Sigma |\vec{H}|^2\, d\mathcal{H}^n .
\end{equation}

Taking limits in the flux identity therefore yields
\begin{equation} 
\int_\Sigma |\vec{H}|^2\, d\mathcal{H}^n = 0 .
\end{equation}

Since $|\vec{H}|^2$ is nonnegative and continuous, this implies $|\vec{H}|^2 \equiv 0$ on $\Sigma$. Hence $\vec{H} \equiv 0$ and $\Sigma$ is minimal. The translator equation $\vec{H}=e_{n+1}^\perp$ therefore implies $e_{n+1}^\perp \equiv 0$ along $\Sigma$.

For each end parameterization $F_k(y)=u_k(y)e_1+y$, the conditions $H\equiv0$, $e_{n+1}^\perp=0$, and $\lim_{|y|\to\infty}\|Du_k(y)\|=0$ imply that $u_k$ is constant by Theorem~\ref{thm:Liouville Rigidity for Asymptotic Minimal Graphs}. Thus each end $E_k$ lies in a vertical hyperplane.

Minimal hypersurfaces are real analytic. Since $\Sigma$ is connected and an open subset of it lies in a vertical hyperplane, analytic continuation implies that the entire hypersurface coincides with that hyperplane.Therefore $\Sigma$ is a vertical hyperplane.
\end{proof}

\section{Conclusion}

In this paper, we have proven that any complete, properly immersed, connected translator in $\mathbb{R}^{n+1}$ ($n \ge 3$) with finite topology and ends $C^1$-asymptotic to vertical hyperplanes must globally be a vertical hyperplane. By evaluating the universal flux identity over a cylindrical exhaustion, we demonstrated that the total squared mean curvature of the surface equals the boundary flux of the translation vector. We showed that this flux vanishes asymptotically by partitioning the parameterizing hyperplane into a canonical wake set and an exponential region to extract strict, bifurcated gradient decay bounds. Consequently, the mean curvature identically vanishes, and the surface reduces to a flat minimal graph via Liouville rigidity.

\noindent\textbf{Mathematics Subject Classification (2020).}

Primary: 53E10 (Mean curvature flow); 35J60 (Nonlinear elliptic equations).

Secondary: 35B40 (Asymptotic behavior of solutions); 35B65 (Smoothness and regularity of solutions)

\medskip
\noindent\textbf{Keywords.}
Mean curvature flow; translating solitons; translators; elliptic partial differential equations; drift Laplacian; anisotropic decay; asymptotic analysis; minimal hypersurfaces; geometric analysis.

\end{document}